\pdfoutput=1
\documentclass[11pt]{article}
\usepackage[cp1251]{inputenc}
\usepackage[english,russian]{babel}
\usepackage{amsmath}
\usepackage{amsthm}
\usepackage{amssymb}
\usepackage{hyperref}
\hypersetup{pdfauthor={K. V. Amelkin, A. V. Kostin},
pdftitle={On asymptotic expansions of oscillating solutions of quasilinear ordinary
    differential equations systems}}

\newtheorem{theorem}{Теорема}
\newtheorem{lemma}{Лемма}
\newtheorem{example}{Пример}
\newtheorem{remark}{Замечание}
\newtheorem{definition}{Определение}

\newcommand{\EditionDate}{24.06.2013}
\allowdisplaybreaks
\date{
Draft of \EditionDate
\\
Черновик от \EditionDate
}

\begin{document}

\title{
On asymptotic expansions of oscillating solutions of quasilinear ordinary
differential equations systems
\\
Об асимптотических разложениях колеблющихся решений квазилинейных систем
обыкновенных дифференциальных уравнений
}
\author{
Kirill Vadimovich Amelkin, Alexander Vasilevich Kostin
\\
Амелькин Кирилл Вадимович, Костин Александр Васильевич
}

\maketitle

\renewcommand{\abstractname}{Abstract}

\begin{abstract}
The existence of a formal particular solution (family of solutions) of the
type \eqref{FSolution} under certain conditions has been proved for the
system of differential equations \eqref{DEV:main}.
The asymptotic nature of this solution (the family of solutions) is
investigated in two individual cases when all the eigenvalues of the matrix
$ A $ 1) are not pure imaginary, 2) are simple and under some additional
assumptions.

-----

Для системы дифференциальных уравнений \eqref{DEV:main} при определённых
условиях установлено существование формального частного решения (семейства
решений) типа \eqref{FSolution}.
Асимптотический характер этого решения (семейства решений) исследован в двух
отдельных случаях, когда все собственные значения матрицы $ A $ 1) не чисто
мнимые, 2) простые и при некоторых дополнительных предположениях.
\end{abstract}


\textbf{Ключевые слова:} квазилинейные системы ОДУ, многочастотные системы
\\

\textbf{UDC (УДК)} 517.928.

\section{Введение}

\nocite{Bib:Kamenec-Podolsky96,Bib:Mastersthesis,Bib:Chernovci98,
Bib:Odessa00,Bib:Chernovci01,Bib:Chernovci02,Bib:Odessa01,Bib:RNASU04}
\nocite{Bib:Kostin_DE65,Bib:Kostin_DE67_I,Bib:Kostin_DE67_II,Bib:Kostin_PHDT,
Bib:Kostin_DE87,Bib:Kostin_DNANU95,Bib:Abgaryan2008,Bib:Mitropolskiy,
Bib:Vazov,Bib:Demidovich,Bib:Malkin,Bib:Levitan,Bib:Sokolov,Bib:Persidskiy,
Bib:Cotton}
\bibliographystyle{gost2008}

В классических трудах Пуанкаре и Мальмквиста (смотрите \cite{Bib:Vazov})
детально исследованы асимптотические разложения решений дифференциальных
уравнений в форме степенных рядов по степеням функции $ \dfrac 1 t $.

В этой работе получены достаточные условия (теорема \ref{T:FSolution})
существования у вещественного векторного квазилинейного обыкновенного
дифференциального уравнения \eqref{DEV:main}, с исчезающими при
$ t \to + \infty $ параметрами $ \varepsilon_l(t) $
$ \left( l = \overline{ 1, m } \right) $ при нелинейностях, формального
частного решения (семейства решений) типа \eqref{FSolution}, которое имеет
вид асимптотического ряда по степеням функций $ \varepsilon_l(t) $
$ \left( l = \overline{ 1, m } \right) $ и их производных любого порядка с
ограниченными на бесконечном промежутке $ I $ коэффициентами, являющимися
векторами из некоторого абстрактного класса
$ \mathcal K_1^{ n \times 1 }(A) $ колеблющихся функций.
Исследован асимптотический характер построенного рекуррентным образом
формального частного решения (семейства решений) \eqref{FSolution} в двух
отдельных случаях, когда все собственные значения матрицы $ A $ 1) не чисто
мнимые (теорема \ref{T:NPI_AsymptChart_O}), 2) простые (теоремы
\ref{T:AsymptChart_O} -- \ref{T:AsymptChart}) и при некоторых дополнительных
предположениях.
В теоремах об асимптотическом характере указано количество произвольных (либо
достаточно малых по абсолютной величине произвольных) скалярных постоянных,
от которых зависит погрешность $ r $.

Новизной данного направления является наличие колеблющихся функций из
некоторого абстрактного класса $ \mathcal K_1(A) $ в векторном
дифференциальном уравнении \eqref{DEV:main} и в асимптотическом разложении
\eqref{FSolution} его формального частного решения (семейства решений).
Публикации других математиков о рассматриваемом здесь случае авторам не
встречались.

Это исследование примыкает также к известным результатам метода малого
параметра Пуанкаре (смотрите \cite{Bib:Malkin}).

\section{Обозначения}

Введем некоторые обозначения, которые будут использоваться всюду в этой
работе.
Будем считать, что множество
$ \mathbb N_0 := \{ 0 \} \cup \mathbb N $,
$ \mathbb R^- := ( - \infty, 0 ) $,
$ \mathbb R_0^- := \mathbb R^- \cup \{ 0 \} $,
$ \mathbb R^+ := ( 0, + \infty ) $,
$ \mathbb R_0^+ := \{ 0 \} \cup \mathbb R^+ $,
промежуток $ I := [ t_0, + \infty ) \subset \mathbb R $,
$ [ a ] $ и $ \{ a \} $ --- целая и дробная части числа $ a \in \mathbb R $
соответственно, $ \Re \, z $ и $ \Im \, z $ --- вещественная и мнимая части числа
$ z \in \mathbb C $ соответственно, $ \imath $ --- мнимая единица, $ E $ ---
единичная матрица, $ \mathbb C^0 \equiv \emptyset $.
Пусть функция $ f(t) : I \to \mathbb R $, число $ a \in \mathbb R $, тогда
символы $ f(t) \downarrow \uparrow a $ ($ t \in I $, $ t \to + \infty $)
будут указывать нам, что $ f(t) $, по крайней мере, нестрого монотонна на
промежутке $ I $ и стремится к $ a $ при аргументе $ t \to + \infty $.
Пусть матрица $ A \in \mathbb C^{ n \times n } $, тогда условимся обозначать
через $ \lambda_j(A) $ $ \left( j = \overline{ 1, n } \right) $ её
собственные значения, $ \Lambda( A ) $ --- её спектр, а множество
$ \Delta( A ) := \left\{ \left. \lambda_j(A) - \lambda_l(A) \, \right|
j, l = \overline{ 1, n }, \, j \ne l \right\} $ (порядок $ n \ge 2 $).
Пусть вектор $ a \in \mathbb C^{ n \times 1 }, $ матрица
$ A \in \mathbb C^{ n \times m }, $ тогда
$ \left( a \right)_j $ --- $ j $
компонента $ a \ \left( j = \overline{ 1, n } \right), $
$ ( A )_{ j l } $ --- элемент $ A $ с индексами
$ j , \, l \left( j = \overline{ 1, n }, \ l = \overline{ 1, m } \right), $
$ A^T $ --- транспонированная $ A $.
Пусть числа $ a_1, a_2, \dots, a_n \in \mathbb C $, тогда символами
$ \operatorname{ diag } ( a_1, \dots, a_n ) $ будем обозначать
диагональную матрицу на главной диагонали, которой стоят элементы
$ a_1, a_2, \dots, a_n $.
Пусть матрица
$ A( x ) : \mathbb R^{ k \times 1 } \to \mathbb C^{ n \times m } $,
$ S $ --- некоторый класс скалярных функций, тогда
под записью $ A( { x } ) \in S $ будем понимать, что элементы
$ \left( A \left( { x } \right) \right)_{ j l } \in S $
$ \left( j = \overline{ 1, n }, \ l = \overline{ 1, m } \right). $
Нормой матрицы $ A \in \mathbb C^{ n \times m }, $ будем называть число
\[
|| A || := \max_{ \substack{ j = \overline{ 1, n } \\
l = \overline{ 1, m } } } \left| ( A )_{ j l } \right|.
\]
Пусть число $ k \in \mathbb N $, матрица
$ A( t ) : \mathbb R \to \mathbb C^{ n \times m } $, тогда обозначим
\[
\left( A(t) \right)^{ (k) } := \left( \frac{ d^k \left( \left( A(t)
\right)_{ jl } \right) }{ d t^k } \right)_{ j, \, l = 1 }^{ n, \, m },
\ \left( A(t) \right)' := \left( A(t) \right)^{ (1) }.
\]
Пусть мультииндекс $ \wp \in \mathbb N_0^n $, вектор
$ y \in \mathbb R^{ n \times 1 } $, вектор--функция
$ { f } ( t, { y } ) :
\mathbb R^{ 1+n } \to \mathbb C^{ m \times 1 } $,
функция $ g ( t, { y } ) : \mathbb R^{ 1+n } \to \mathbb C $,
тогда
\begin{gather*}
| \wp | := \sum_{ j=1 }^n ( \wp )_j,
\ \ \wp ! := \prod_{ j=1 }^n ( \wp )_j !,
\ \ { y }^{ \, \wp } := \prod_{ j=1 }^n ( y )_j^{ ( \wp )_j },
\\
\partial_{ { y } }^\wp { f } ( t, { y } ) :=
\frac { \partial^{ | \wp | } { f } ( t, { y } ) }
{ \partial ( y )_1^{ ( \wp )_1 } \partial ( y )_2^{ ( \wp )_2 } \dots
\partial ( y )_n^{ ( \wp )_n } },
\ \nabla_{ y } \, g ( t, { y } ) :=
\left( \frac{ \partial g ( t, { y } ) }{ \partial ( y )_1 }, \dots,
\frac{ \partial g ( t, { y } ) }{ \partial ( y )_n } \right).
\end{gather*}

Везде в этой статье будем полагать, что сумма
$ \sum\limits_{ j = k_0 }^{ k_1 } a_j = 0 $,
произведение
$ \prod\limits_{ j = k_0 }^{ k_1 } a_j = 1 $,
$ \overline{ k_0, k_1 } = \emptyset $,
функция
$ { f } (b_1, \dots, b_{ k_2 }, a_{ k_0 }, \dots, a_{ k_1 },
c_1, \dots, c_{ k_3 }) \equiv
{ f } (b_1, \dots, b_{ k_2 }, c_1, \dots, c_{ k_3 }) $
при $ k_1 < k_0 $, $ f( t, c ) \equiv f(t) $, если $ c \in \emptyset $.

В этой публикации фигурируют производные только натурального порядка.

\section{Существование формального частного решения (семейства решений)}

В этой секции статьи построено, по крайней мере, одно формальное частное
решение вида \eqref{FSolution} основного вещественного векторного
квазилинейного обыкновенного дифференциального уравнения \eqref{DEV:main}
(теорема \ref{T:FSolution}).

\subsection{\texorpdfstring{Абстрактный класс $ \mathcal K_1(A) $
колеблющихся функций}{Первый абстрактный класс колеблющихся функций}}

Пусть, далее, матрица $ A \in \mathbb R^{ n \times n } $, порядок
$ n \ge 2 $.

В работах \cite{Bib:Kamenec-Podolsky96} -- \cite{Bib:Odessa00} авторами
настоящей статьи при некоторых условиях был исследован асимптотический
характер формального частного решения типа \eqref{FSolution} (с
периодическими векторными коэффициентами) векторного дифференциального
уравнения вида \eqref{DEV:main} в случае, когда $ f(t) $,
$ \partial_{ y }^\wp f _l( t, y ) $
$ \left( \wp \in \mathbb N_0^n, \ l = \overline{ 1, m } \right) $ являются
периодическими вектор--функциями по независимой переменной $ t $.
С целью обобщения ранее полученных результатов на более общий вид
вектор--функций $ f(t) $, $ \partial_{ y }^\wp f _l( t, y ) $
$ \left( \wp \in \mathbb N_0^n, \ l = \overline{ 1, m } \right) $ в
исследованиях \cite{Bib:Chernovci01} -- \cite{Bib:RNASU04} авторами данной
публикации были аксиоматически введены в рассмотрение абстрактные классы
$ \mathcal K_1(A) $, $ \mathcal K_2(A) $ и $ \mathcal K_3(A) $ колеблющихся
функций.
Ниже дано аксиоматическое определение класса $ \mathcal K_1(A) $
необходимого, для того чтобы найти векторные коэффициенты формального
частного решения (семейства решений) \eqref{FSolution} векторного
дифференциального уравнения \eqref{DEV:main}.

\begin{definition}

Будем обозначать через $ \mathcal K_1(A) $ некоторый непустой класс функций
$ \left\{ f(t) : I \to \mathbb R \right\} $ таких, что
\begin{enumerate}
\item
$ f(t) \in C( I ) $, супремум
$ \sup\limits_{ t \in I } \left| f(t) \right| < + \infty $;
                                                            \label{C:1_aca}
\item
если функции $ f_1(t), f_2(t) \in \mathcal K_1(A) $, $ \alpha \in \mathbb R $,
то $ \alpha f_1(t) $, $ f_1(t) + f_2(t) $, $ f_1(t) f_2(t) $
$ \in \mathcal K_1(A) $;
                                                            \label{C:2_aca}
\item
пусть функция $ \tilde f(t) : I \to \mathbb C $ и $ \Re \, \tilde f(t) $,
$ \Im \, \tilde f(t) \in \mathcal K_1(A) $, тогда у каждого дифференциального
уравнения вида
\begin{equation}
\frac{ dy }{ dt } = \lambda y + \tilde f(t)
\ \ ( y : I \to \mathbb C, \ t \in I, \ \lambda \in \Lambda( A ) )
                                        \label{DE:splitted_truncated_main}
\end{equation}
на промежутке $ I $ существует, по крайней мере, одно частное решение
$ y( t, \tilde c ) $ такое, что $ \Re \, y( t, \tilde c ), $
$ \Im \, y( t, \tilde c ) \in \mathcal K_1(A) \cap C^1(I) $, параметр
$ \tilde c \in \mathbb C^{ \tilde n } $ ($ 0 \le \tilde n \le 1 $).
                                                            \label{C:3_aca}
\end{enumerate}
                                                            \label{D:aca}
\end{definition}

Приведем три важных примера классов $ \mathcal K_1(A) $ состоящих из
\begin{enumerate}
\item
периодических функций,
\item
равномерных почти периодических в смысле Г. Бора (далее в тексте этой работы
РПП) функций с конечными спектрами,
\item
конечных линейных комбинаций функций вида $ e^{ \gamma t } $ (число
$ \gamma \in \mathbb R_0^- ( \mathbb R^- ) $) с действительными
коэффициентами.
\end{enumerate}

\begin{example}
                                                    \label{Ex:PeriodicK1}
Пусть число $ \tau \in \mathbb R^+ $ фиксировано, множество
\begin{equation}
\Gamma_1 := \left\{ \left. \dfrac{ 2 \pi \imath s }{ \tau } \, \right|
s \in \mathbb Z \right\},
\ \Lambda( A ) \cap \Gamma_1 = \emptyset.
                                        \label{C:FourierIndex_EV_separated}
\end{equation}
Тогда множество вещественных непрерывно дифференцируемых на $ \mathbb R $
$ \tau $\nobreakdash-\hspace{0mm}периодических функций образует некоторый класс
$ \mathcal K_1(A) $ на любом промежутке $ I $.
Причём у каждого дифференциального уравнения вида
\eqref{DE:splitted_truncated_main} $ ( \lambda \in \Lambda( A ) ) $ в таком
случае будет существовать ровно по одному частному решению, вещественная и
мнимая части которых будут принадлежать классу
$ \mathcal K_1(A) \cap C^1(I) $.

Очевидно, что для указанного множества функций справедливы аксиомы
\ref{C:1_aca}., \ref{C:2_aca}. класса $ \mathcal K_1(A) $ на любом промежутке
$ I $.
Убедимся, что также для этого множества выполняется и аксиома \ref{C:3_aca}.
класса $ \mathcal K_1(A) $ на любом промежутке $ I $.
Действительно, в силу теоремы о порядке убывания коэффициентов Фурье
периодической функции в зависимости от её гладкости, при указанных условиях
неоднородная часть дифференциального уравнения
\eqref{DE:splitted_truncated_main} разлагается в абсолютно и равномерно
сходящийся на любом промежутке $ I $ ряд Фурье
$$
\tilde f(t) = \sum\limits_{ s = - \infty }^\infty a_s
e^{ \dfrac{ 2 \pi \imath s }{ \tau } t },
\ \ a_s = \frac 1 \tau \int\limits_0^\tau \tilde f(t)
e^{ - \dfrac { 2 \pi \imath s } \tau t } \, d t \ \ (s \in \mathbb Z),
\ \ \sum\limits_{ s = - \infty }^\infty | a_s | < + \infty.
$$
В таком случае непосредственной подстановкой ряда Фурье
$$
y( t ) := \sum\limits_{ s = - \infty }^\infty
\dfrac{ a_s }{ \dfrac{ 2 \pi \imath s }{ \tau } - \lambda }
e^{ \dfrac{ 2 \pi \imath s }{ \tau } t } \ \ (\lambda \in \Lambda( A ))
$$
в дифференциальное уравнение \eqref{DE:splitted_truncated_main} можно
убедиться, что он является частным решением этого уравнения на любом
промежутке $ I $.
Очевидно, что функции
$ \Re \, y( t ), $ $ \Im \, y( t ) \in \mathcal K_1(A) \cap C^1(I) $.
Заметим, что при условиях этого примера общим решением дифференциального
уравнения \eqref{DE:splitted_truncated_main} будет функция
\[
y( t, \tilde c ) := \tilde c e^{ \lambda t } + y( t )
\ \ ( \tilde c \in \mathbb C, \ \lambda \in \Lambda( A ), \ t \in I ).
\]
Однако при $ \tilde c \ne 0 $ ($ \forall \, \lambda \in \mathbb C $) функции
$ \Re \, y( t, \tilde c ), \ \Im \, y( t, \tilde c ) $ не могут принадлежать классу
$ \mathcal K_1(A) $, так как при этом нарушаются условия
$ \tau $\nobreakdash-\hspace{0mm}периодичности,
\eqref{C:FourierIndex_EV_separated}.

\end{example}

Примеру РПП класса $ \mathcal K_1(A) $ предпошлем пример простейшей РПП
функции.

\begin{example}

Функция $ \sin( \gamma_1 t ) + \sin( \gamma_2 t ) $
$(\gamma_1, \gamma_2, t \in \mathbb R, \ \gamma_1 \gamma_2 \ne 0) $ является
(чисто) периодической, если $ \dfrac{ \gamma_1 }{ \gamma_2 } \in \mathbb Q $
и РПП (непериодической), если
$ \dfrac{ \gamma_1 }{ \gamma_2 } \in \mathbb R \setminus \mathbb Q $.

\end{example}

Ниже используем известные понятия (смотрите \cite{Bib:Demidovich},
\cite{Bib:Levitan}).

\begin{definition}

Средним значением функциональной матрицы
$ F(t) : I \to \mathbb C^{ m \times n } $ называется матрица
\[
\operatorname{ M } ( F(t) ) := \lim_{ t \to + \infty }
\frac 1 t \int\limits_{ t_0 }^t F( \tau ) \, d \tau,
\]
если $ F(t) $ интегрируема на любом отрезке $ [ t_0, t_1 ] $
$ ( t_1 > t_0 ), $ и если указанный здесь предел существует и конечен.

\end{definition}

\begin{definition}

Рядом Фурье РПП матрицы $ F(t) : \mathbb R \to \mathbb C^{ m \times n } $
называется конечный или счётный тригонометрический матричный ряд
\[
F(t) \sim \sum_s A_s e^{ \imath \gamma_s t },
\ \ A_s = \operatorname{ M } \left( F(t) e^{ - \imath \gamma_s t } \right),
\]
где матрицы $ A_s $ --- коэффициенты Фурье, $ \gamma_s $ --- показатели Фурье,
множество $ \left\{ \gamma_s \right\} \subset \mathbb R $ --- спектр
матричной функции $ F(t) $.

\end{definition}

\begin{example}
                                                \label{Ex:AlmostPeriodicK1}
Пусть задана конечная или счётная совокупность действительных чисел
$ \Omega_1 := \left\{ \omega_1, \ \omega_2, \dots \right\} \subset
\mathbb R $,
множество $ \Gamma_2 $ состоит из всевозможных конечных линейных комбинаций
чисел из $ \Omega_1 $ с целыми коэффициентами, то есть
\begin{equation*}
\Gamma_2 := \left\{ \left. \sum_{ s=1 }^p k_s \omega_s \right| k_s \in
\mathbb Z,
\ \omega_s \in \Omega_1, \ s = \overline{ 1, p }, \ p \in \mathbb N \right\},
                                                            \label{Gamma_2}
\end{equation*}
причём выполнено условие
\begin{equation}
\inf_{ \substack{ \gamma \in \Gamma_2 \\ \lambda \in \Lambda( A ) } }
\left| \imath \gamma - \lambda \right| > 0.
                                            \label{C:spectrum_EV_separated}
\end{equation}
Тогда множество вещественных РПП функций обладающих конечными спектрами,
являющимися подмножествами $ \Gamma_2 $, то есть
\begin{equation}
\left\{ \left. \sum\limits_{ s = 1 }^p a_s e^{ \imath \gamma_s t } +
\bar a_s e^{ - \imath \gamma_s t } \right|
a_s \in \mathbb C, \ \gamma_s \in \Gamma_2,
\ s = \overline{ 1, p }, \ p \in \mathbb N \right\}  \ \ ( t \in I ),
                                                \label{Set:AlmostPeriodicK_1}
\end{equation}
образует некоторый класс $ \mathcal K_1(A) $ на любом промежутке $ I $.
Причём у каждого дифференциального уравнения вида
\eqref{DE:splitted_truncated_main} $ ( \lambda \in \Lambda( A ) ) $ в таком
случае будет существовать ровно по одному частному решению, вещественная и
мнимая части которых будут принадлежать классу
$ \mathcal K_1(A) \cap C^1(I) $.

Очевидно, что для множества \eqref{Set:AlmostPeriodicK_1} справедливы аксиомы
\ref{C:1_aca}., \ref{C:2_aca}. класса $ \mathcal K_1(A) $ на любом промежутке
$ I $.
Убедимся, что также для этого множества выполняется и аксиома \ref{C:3_aca}.
класса $ \mathcal K_1(A) $ на любом промежутке $ I $.
Действительно, при указанных условиях неоднородную часть дифференциального
уравнения \eqref{DE:splitted_truncated_main} можно представить в виде
$$
\tilde f(t) = \sum\limits_{ s = 1 }^p a_s e^{ \imath \gamma_s t } +
b_s e^{ - \imath \gamma_s t }
\ \ \left( a_s, \ b_s \in \mathbb C, \ \gamma_s \in \Gamma_2,
\ s = \overline{ 1, p }, \ t \in I \right).
$$
В таком случае непосредственной подстановкой функции
$$
y( t ) := \sum\limits_{ s = 1 }^p
\frac{ a_s }{ \imath \gamma_s - \lambda } e^{ \imath \gamma_s t } -
\frac{ b_s }{ \imath \gamma_s + \lambda } e^{ - \imath \gamma_s t }
\ \ (\lambda \in \Lambda( A ))
$$
в дифференциальное уравнение \eqref{DE:splitted_truncated_main} можно
убедиться, что она является частным решением этого уравнения на любом
промежутке $ I $.
Очевидно, что функции
$ \Re \, y( t ), $ $ \Im \, y( t ) \in \mathcal K_1(A) \cap C^1(I) $.
Заметим, что при условиях этого примера общим решением дифференциального
уравнения \eqref{DE:splitted_truncated_main} будет функция
\[
y( t, \tilde c ) := \tilde c e^{ \lambda t } + y( t )
\ \ ( \tilde c \in \mathbb C, \ \lambda \in \Lambda( A ), \ t \in I ).
\]
Однако при $ \tilde c \ne 0 $ ($ \forall \, \lambda \in \mathbb C $) функции
$ \Re \, y( t, \tilde c ), \ \Im \, y( t, \tilde c ) $ не могут принадлежать классу
$ \mathcal K_1(A) $, так как при этом нарушаются условия
\eqref{Set:AlmostPeriodicK_1}, \eqref{C:spectrum_EV_separated}.

\end{example}

Если в предыдущих примерах классов $ \mathcal K_1(A) $ функции колебались на
всём множестве $ \mathbb R $ около своего среднего значения, то в
нижеследующем они монотонно (для достаточно больших значений аргумента $ t $)
стремятся к нему (при $ t \to + \infty $).

\begin{example}
                                            \label{Ex:Decreasing_Exp_LC_K1}
Пусть задана конечная или счётная совокупность неположительных (либо
отрицательных) действительных чисел
$ \Omega_2 := \left\{ \omega_1, \ \omega_2, \dots \right\} \subset
\mathbb R_0^- \ ( \mathbb R^- ) $,
множество $ \Gamma_3 $ состоит из всевозможных конечных нетривиальных
линейных комбинаций чисел из $ \Omega_2 $ с коэффициентами из множества
$ \mathbb N_0 $, то есть
\begin{equation*}
\Gamma_3 := \left\{ \left. \sum_{ s=1 }^p k_s \omega_s \right|
k_s \in \mathbb N_0, \ \sum_{ s=1 }^p k_s \ne 0,
\ \omega_s \in \Omega_2,
\ s = \overline{ 1, p }, \ p \in \mathbb N \right\},
\end{equation*}
причём выполнено условие
\begin{equation}
\inf_{ \substack{ \gamma \in \Gamma_3 \\ \lambda \in \Lambda( A ) } }
\left| \gamma - \lambda \right| > 0.
                                            \label{C:exponent_EV_separated}
\end{equation}
Тогда множество конечных линейных комбинаций
\begin{equation}
\left\{ \left. \sum_{ s=1 }^p a_s e^{ \gamma_s t } \right|
a_s \in \mathbb R, \ \gamma_s \in \Gamma_3,
\ s = \overline{ 1, p }, \ p \in \mathbb N \right\} \ \ ( t \in I )
                                                \label{Set:Decreasing_Exp_LC}
\end{equation}
образует некоторый класс $ \mathcal K_1(A) $ на любом промежутке $ I $
(константа $ t_0 \ge 0 $).
Причём у каждого дифференциального уравнения вида
\eqref{DE:splitted_truncated_main} $ ( \lambda \in \Lambda( A ) ) $ в таком
случае будет существовать ровно по одному частному решению, вещественная и
мнимая части которых будут принадлежать классу
$ \mathcal K_1(A) \cap C^1(I) $.

Очевидно, что для множества \eqref{Set:Decreasing_Exp_LC} справедливы аксиомы
\ref{C:1_aca}., \ref{C:2_aca}. класса $ \mathcal K_1(A) $ на любом промежутке
$ I $ (константа $ t_0 \ge 0 $).
Убедимся, что также для этого множества выполняется и аксиома \ref{C:3_aca}.
класса $ \mathcal K_1(A) $ на любом промежутке $ I $ (константа
$ t_0 \ge 0 $).
Действительно, при указанных условиях неоднородную часть дифференциального
уравнения \eqref{DE:splitted_truncated_main} можно представить в виде
$$
\tilde f(t) = \sum\limits_{ s=1 }^p a_s e^{ \gamma_s t }
\ \ \left( a_s \in \mathbb C, \ \gamma_s \in \Gamma_3,
\ s = \overline{ 1, p }, \ t \in I \right).
$$
В таком случае непосредственной подстановкой функции
\[
y( t ) =: \sum\limits_{ s=1 }^p
\frac{ a_s }{ \gamma_s - \lambda } e^{ \gamma_s t }
\ \ ( \lambda \in \Lambda( A ) )
\]
в дифференциальное уравнение \eqref{DE:splitted_truncated_main} можно
убедиться, что она является частным решением этого уравнения на любом
промежутке $ I $ (константа $ t_0 \ge 0 $).
Очевидно, что функции
$ \Re \, y( t ), \ \Im \, y( t ) \in \mathcal K_1(A) \cap C^1(I) $.
Заметим, что при условиях этого примера общим решением дифференциального
уравнения \eqref{DE:splitted_truncated_main} будет функция
\[
y( t, \tilde c ) := \tilde c e^{ \lambda t } + y( t )
\ \ ( \tilde c \in \mathbb C, \ \lambda \in \Lambda( A ), \ t \in I ).
\]
Однако при $ \tilde c \ne 0 $ ($ \forall \, \lambda \in \mathbb C $) функции
$ \Re \, y( t, \tilde c ), \ \Im \, y( t, \tilde c ) $ не могут принадлежать классу
$ \mathcal K_1(A) $, так как при этом нарушаются условия
\eqref{Set:Decreasing_Exp_LC}, \eqref{C:exponent_EV_separated}.

\end{example}

В следующей лемме доказано существование, по крайней мере, одного частного
решения из класса $ \mathcal K_1^{ n \times 1 }(A) \cap C^1( I ) $
укороченного уравнения \eqref{DEV:truncated_main} соответствующего основному
векторному обыкновенному дифференциальному уравнению \eqref{DEV:main}.

\begin{lemma}

Пусть в векторном дифференциальном уравнении
                                                    \label{L:truncated_main}
\begin{equation}
\frac{ d \varphi }{ dt } = A { \varphi } + { f }(t) \ \ ( t \in I )
                                                   \label{DEV:truncated_main}
\end{equation}
\begin{enumerate}
\item
вектор--функция $ \varphi : I \to \mathbb R^{ n \times 1 } $, порядок
$ n \ge 2 $;
\item
матрица $ A \in \mathbb R^{ n \times n } $ и
задан некоторый класс $ \mathcal K_1(A) $;
\item
вектор--функция $ f(t) \in \mathcal K_1^{ n \times 1 }(A). $
\end{enumerate}
Тогда у векторного дифференциального уравнения \eqref{DEV:truncated_main}
на промежутке $ I $ существует, по крайней мере, одно частное решение
$ \varphi_0( t, c_0 ) $ из класса
$ \mathcal K_1^{ n \times 1 }(A) \cap C^1( I ) $ (параметр
$ c_0 \in \mathbb R^{ n_0 } $, $ 0 \le n_0 \le n $).

\end{lemma}

\begin{proof}
Приведем матрицу $ A $ к жордановой канонической форме $ J $.
Для этого сделаем в векторном дифференциальном уравнении
\eqref{DEV:truncated_main} замену неизвестной вектор--функции вида
$ \varphi = P_0 z, $ где матрица $ P_0 \in \mathbb C^{ n \times n } $,
$ J = P_0^{ -1 } A P_0 $, определитель $ \det P_0 \ne 0 $.
В результате получим векторное дифференциальное уравнение вида
\begin{equation}
\frac{ dz }{ dt } = J  z + P_0^{ -1 }  f(t).
                                \label{DEV:JNF_trunc_main}
\end{equation}

Учитывая структуру жордановой канонической формы матрицы, заметим, что решать
систему $ n $ скалярных дифференциальных уравнений \eqref{DEV:JNF_trunc_main}
нужно снизу вверх.
При такой последовательности решения дифференциальной системы
\eqref{DEV:JNF_trunc_main} каждый раз будем интегрировать одно
скалярное дифференциальное уравнение вида
$$
\frac{ d \tilde z }{ dt } = \lambda \tilde z + \tilde f(t)
\ \ \left( \lambda \in \Lambda( A ) \right)
$$
с уже известным из предыдущих шагов свободным членом $ \tilde f(t) $.

В силу аксиом класса $ \mathcal K_1(A) $ вектор--функции
$ \Re \, P_0^{ -1 } f(t), $
$ \Im \, P_0^{ -1 } f(t) \in \mathcal K_1^{ n \times 1 }(A). $
Поэтому, в силу этих же аксиом, на промежутке $ I $ существует, по крайней
мере, одно частное решение $  z_0( t, \tilde c_0 ) $ (параметр
$ \tilde c_0 \in \mathbb C^{ n_0 } $, $ 0 \le n_0 \le n $) векторного
дифференциального уравнения \eqref{DEV:JNF_trunc_main}, такое, что
вектор--функции $ \Re \, z_0( t, \tilde c_0 ) $,
$ \Im \, z_0( t, \tilde c_0 ) \in
\mathcal K_1^{ n \times 1 }(A) \cap C^1(I) $.

Возвращаясь обратно к неизвестной вектор--функции $  \varphi $, получим,
по крайней мере, одно комплекснозначное частное решение
$ P_0  z_0( t, \tilde c_0 ) $
векторного дифференциального уравнения \eqref{DEV:truncated_main}, такое, что
$ \Re \, P_0 z_0( t, \tilde c_0 ), $
$ \Im \, P_0  z_0( t, \tilde c_0 ) \in
\mathcal K_1^{ n \times 1 }(A) \cap C^1(I). $
Записывая вектор-функцию $ P_0 z_0( t, \tilde c_0 ) $ в виде выражения
$ \Re \, P_0 z_0( t, \tilde c_0 ) + \imath \Im \, P_0 z_0( t, \tilde c_0 ) $
и подставляя его в исходное векторное дифференциальное уравнение
\eqref{DEV:truncated_main}, убедимся, что
$ \varphi_0( t, c_0 ) := \Re \, P_0 z_0( t, \tilde c_0 ) $
является его частным решением (семейством решений) из класса
$ \mathcal K_1^{ n \times 1 }(A) \cap C^1( I ) $ (параметр
$ c_0 \in \mathbb R^{ n_0 } $).

\end{proof}

\subsection{Основное вещественное векторное квазилинейное обыкновенное
дифференциальное уравнение}

Основным объектом исследований в настоящей работе является вещественное
векторное квазилинейное обыкновенное дифференциальное уравнение
\begin{equation}
\frac{ dy }{ dt } = A { y } + { f }(t) +
\sum_{ l=1 }^m \varepsilon_l(t) { f }_l( t, { y } )
\ \ \left( ( t, { y } ) \in D \right),
                                                        \label{DEV:main}
\end{equation}
в предположении, что
\begin{enumerate}
\item
выполнены все условия леммы \ref{L:truncated_main};
                                                    \label{C:main_t_first}
\item
область
$ D := \left\{ ( t, { y } ) \ \left| \ t \in I,
\ { y } \in \mathbb R^{ n \times 1 }, \right.
|| { y } - \varphi_0( t, c_0 ) || \le a \right\} $
(число $ a \in \mathbb R^+ $), причём вектор--функция $ \varphi_0( t, c_0 ) $
(параметр $ c_0 \in \mathbb R^{ n_0 } $, $ 0 \le n_0 \le n $) --- некоторое
частное решение из класса $ \mathcal K_1^{ n \times 1 }(A) \cap C^1( I ) $
укороченного векторного дифференциального уравнения
\eqref{DEV:truncated_main} соответствующего основному уравнению
\eqref{DEV:main};
                                                    \label{C:main_t_second}
\item
функции $ \varepsilon_l(t) : I \to \mathbb R $,
$ \varepsilon_l(t) \in C^\infty ( I ) $,
$ \dfrac{ d^s \varepsilon_l(t) }{ d t^s } = o(1) $ ($ t \to + \infty $,
$ l = \overline{ 1, m } $, $ s \in \mathbb N_0 $);
\item
вектор--функции
$ { f }_l( t, { y } ) : D \to \mathbb R^{ n \times 1 } $,
$ { f }_l( t, { y } ) \in
C^{ 0, \ \infty }_{ t, \ { y } }( D ) $ ($ l = \overline{ 1, m } $);
\item
супремумы
$ \sup\limits_{ ( t,  y ) \in D } \left\|
\partial_{ { y } }^\wp { f }_l( t, { y } ) \right\| < + \infty, $
$ \partial_{ { y } }^\wp { f }_l( t, \varphi_0( t, c_0 ) )
\in \mathcal K_1^{ n \times 1 }(A) $
$ \left( \wp \in \mathbb N_0^n, \ l = \overline{ 1, m } \right). $
                                                        \label{C:main_t_last}
\end{enumerate}

\subsection{Ранг функций}

В работах \cite{Bib:Kamenec-Podolsky96} -- \cite{Bib:Mastersthesis} авторами
настоящей публикации при некоторых условиях было построено формальное частное
решение одного квазилинейного дифференциального уравнения второго порядка,
которое является частным случаем векторного дифференциального уравнения
\eqref{DEV:main}.
В качестве функций $ \varepsilon_l(t) $
$ \left( l = \overline{ 1, m } \right) $ там рассматривалась одна функция
$ \dfrac 1 t $.
С целью обобщения полученного результата на более общий вид функций в
исследованиях \cite{Bib:Chernovci98} -- \cite{Bib:RNASU04} авторами данной
работы было использовано понятие ранга.
Ниже дано аксиоматическое определение понятие ранга, которое необходимо для
частичного упорядочивания, нумерации с помощью только двух натуральных
индексов и относительной оценки скорости стремления к нулю (в виде
$ O \left( \left( \varepsilon_1(t) \right)^\alpha \right) $,
$ \alpha \in \mathbb R^+ $, $ t \to + \infty $) всевозможных произведений
функций $ \varepsilon_l(t) $ $ \left( l = \overline{ 1, m } \right) $ и их
производных любого натурального порядка.

Каждой из функций $ \varepsilon_l(t) $
$ \left( l = \overline{ 1, m } \right) $ будем приписывать некоторый ранг
$ \varrho_l \in \mathbb R^+ $ $ \left( l = \overline{ 1, m } \right) $,
$ \varrho_{ l-1 } \le \varrho_l $
$ \left( l = \overline{ 2, m }, \ m \ge 2 \right). $
Условимся писать
\[
\operatorname{ Rank } ( \varepsilon_l(t) ) :=
\varrho_l \ \ \left( l = \overline{ 1, m } \right).
\]

\begin{example}
                                                            \label{E:rang}
В случае, когда
\[
\varepsilon_l(t) = \frac { \ln^{ \alpha_{ l1 } } t
\ln^{ \alpha_{ l2 } } \ln t \dots
\ln^{ \alpha_{ l { p_l } } } ( \ln ( \dots ( \ln t ) \dots
) ) } { t^{ \alpha_{ l0 } } }
\ \ (l = \overline{ 1, m }),
\]
где $ \alpha_{ lq } \in \mathbb R $, $ \alpha_{ l0 } > 0 $
($ l = \overline{ 1, m } $, $ q = \overline{ 0, p_l } $)
($ t_0 $ достаточно большое число, чтобы последние дроби имели смысл) удобно
положить, что
\[
\operatorname{ Rank } ( \varepsilon_l(t) ) = \alpha_{ l0 }
\ \ (l = \overline{ 1, m }).
\]

\end{example}

\begin{remark}

Очевидно, что этим примером функций $ \varepsilon_l(t) $
($ l = \overline{ 1, m } $) не исчерпываются все возможные их варианты.
Действительно, при некоторых условиях можно добавить в них в качестве
сомножителей функции вида $ \sin \left( \varepsilon_l(t) \right), $
$ \cos \left( \varepsilon_l(t) \right), $
$ \exp \left( t^\alpha \right) $ ($ l \in \overline{ 1, m } $,
$ \alpha \in \mathbb R^- $).
Тем самым применённое в работе понятие ранга позволяет охватить достаточно
широкий класс функций.

\end{remark}

Определим ранг произведений и производных функций $ \varepsilon_l(t) $
($ l = \overline{ 1, m } $) следующим образом.

\begin{definition}
                                                    \label{D:ProdDerivRank}
Пусть заданы числа $ k_l \in \mathbb N_0 $, $ \beta_{ lr } \in \mathbb N_0 $
($ l = \overline{ 1, m } $, $ r = \overline{ 0, k_l } $),
причём хотя бы одно из чисел $ \beta_{ lr } $ отлично от
нуля и функция
\begin{equation*}
\nu_\eta(t) := \prod_{ l=1 }^m \prod_{ r=0 }^{ k_l }
\left( \dfrac{ d^r \varepsilon_l(t) }{ d t^r } \right)^{ \beta_{ lr } },
\ \eta = ( k_1, \dots, k_m,
\beta_{ 10 }, \dots, \beta_{ m k_m } ),
\end{equation*}
тогда условимся считать, что
\begin{equation}
\operatorname{ Rank } \left( \nu_\eta(t) \right) :=
\sum_{ l=1 }^m \sum_{ r=0 }^{ k_l }
\left( \varrho_l + r \right) \beta_{ lr }.
                                                        \label{ProdDerivRank}
\end{equation}

\end{definition}

Множество всех введённых таким образом рангов обозначим
через $ R $.
\[
R := \left\{ \sum_{ l=1 }^m \sum_{ r=0 }^{ k_l }
\left( \varrho_l + r \right) \beta_{ lr } : k_l,
\beta_{ lr } \in \mathbb N_0 \ (l = \overline{ 1, m },
\ r = \overline{ 0, k_l }),
\sum_{ l=1 }^m \sum_{ r=0 }^{ k_l } \beta_{ lr } > 0 \right\}
\]
Это множество, когда $ \varrho_l $ ($ l = \overline{ 1, m } $) и $ m $
фиксированы, в любом ограниченном интервале
$ (b, c) \subset \mathbb R $ содержит конечное число элементов.
Таким образом множество $ R $ счётно.
Следовательно, можем перенумеровать его
\[
R = \{ \rho_s, \ \rho_s < \rho_{ s+1 } \ ( s \in \mathbb N ) \}.
\]
Заметим, также, что множество $ R $ содержит наименьший
элемент --- $ \varrho_1 $.

\subsection{Теорема о существовании формального частного решения (семейства
решений)}

Для каждого фиксированного $ s \in \mathbb N $ рассмотрим множество
\[
M_s := \left\{ \nu_\eta(t) \bigm|
\operatorname{ Rank } ( \nu_\eta(t) ) = \rho_s \right\}.
\]
Перенумеруем элементы этого множества с помощью индекса
$ p $ и запишем $ M_s $ в форме
\[
M_s = \left\{ \nu_{ sp }(t),
\ p = \overline{ 1, \varkappa_s } \right\} \ ( s \in \mathbb N ),
\]
где $ \varkappa_s $ ($ s \in \mathbb N $) --- количество элементов
множества $ M_s $.

По функциям $ \nu_{ sp }(t) $ ($ s \in \mathbb N $,
$ p = \overline{ 1, \varkappa_s } $) будем производить разложение формального
частного решения (семейства решений) типа \eqref{FSolution} основного
векторного дифференциальное уравнения \eqref{DEV:main}.
Заметим, что их можно представить в виде
\begin{equation}
\nu_{ sp }(t) = \prod_{ l=1 }^m \prod_{ r=0 }^{ [ \rho_s - \varrho_l ] }
\left( \dfrac{ d^r \varepsilon_l(t) }{ d t^r } \right)^{ \beta_{ splr } }
\ ( \beta_{ splr } \in \mathbb N_0 :
\operatorname{ Rank } ( \nu_{ sp }(t) ) = \rho_s ).
                                                                \label{E:nu_sp(t)}
\end{equation}
Очевидно, что функции $ \nu_{ 1p }(t) \equiv \varepsilon_p(t) $
$ \left( p = \overline{ 1, \varkappa_1 } \right) $.

В следующей теореме доказано существование, по крайней мере, одного
формального частного решения вида \eqref{FSolution} у основного векторного
дифференциальное уравнения \eqref{DEV:main}.

\begin{theorem}
                                                        \label{T:FSolution}
Пусть вещественное векторное квазилинейное обыкновенное дифференциальное
уравнение \eqref{DEV:main} удовлетворяет условиям \ref{C:main_t_first}.
-- \ref{C:main_t_last} и заданы некоторые ранги $ \varrho_l $ функций
$ \varepsilon_l(t) $ $ \left( l = \overline{ 1, m } \right) $, для которых
выполняется определение \ref{D:ProdDerivRank}.
Тогда существует, по крайней мере, одно его формальное частное решение вида
\begin{equation}
\tilde y = \varphi_0( t, c_0 ) + \sum_{ s=1 }^\infty
\sum_{ p=1 }^{ \varkappa_s } \nu_{ sp }(t)
\varphi_{ sp }( t, c_{ sp } ) \ \ ( t \in I ),
                                                        \label{FSolution}
\end{equation}
где вектор--функции
$ \varphi_{ sp }( t, c_{ sp } ) \in
\mathcal K_1^{ n \times 1 }(A) \cap C^1( I ) $
(параметры $ c_{ sp } \in \mathbb R^{ n_{ sp } } $, числа
$ n_{ sp } \in \left\{ 0, 1, 2, \dots, n \right\} $,
$ s \in \mathbb N $, $ p = \overline{ 1, \varkappa_s } $) определяются
рекуррентным образом из бесконечной последовательности совокупностей
векторных дифференциальных уравнений вида \eqref{SSDE:Coef}.
(Формальное семейство решений \eqref{FSolution} может зависеть не более чем
от $ n $ произвольных скалярных постоянных.)

\end{theorem}

\begin{proof}

Будем искать формальные решения векторного дифференциального уравнения
\eqref{DEV:main} в виде
\begin{equation}
\tilde y = \varphi_0( t, c_0 ) + \sum_{ s=1 }^\infty
\sum_{ p=1 }^{ \varkappa_s } \nu_{ sp }(t)
\varphi_{ sp }( t ) \ \ ( t \in I ),
                                                        \label{UFSolution}
\end{equation}
где вектор--функции
$ \varphi_{ sp }( t ) \in \mathcal K_1^{ n \times 1 }(A) \cap C^1( I ) $
($ s \in \mathbb N $, $ p = \overline{ 1, \varkappa_s } $) пока не известны.
Формально продифференцируем по независимой переменной $ t $ это равенство, в
результате получим
\begin{equation}
\frac{ d \tilde y }{ dt } =
\frac{ d \varphi_0( t, c_0 ) }{ dt } +
\sum_{ s=1 }^\infty \sum_{ p=1 }^{ \varkappa_s } \left(
\frac{ d \nu_{ sp }(t) }{ dt } \varphi_{ sp }( t ) + \nu_{ sp }(t)
\frac{ d \varphi_{ sp }(t) }{ dt } \right).
                                                        \label{UFSolution'}
\end{equation}

Введем обозначение
$ \tilde \delta(t) := \tilde y - \varphi_0( t, c_0 ). $
С помощью разложения вектор--функций
$ { f }_l ( t, \tilde y ) $ ($ l = \overline{ 1, m } $)
в формальные ряды Тейлора получим
\begin{equation}
\sum_{ l=1 }^m \varepsilon_l(t) { f }_l ( t, \tilde y ) =
\sum \limits _{ l=1 }^m \varepsilon_l(t)
\sum_{ \wp \in \mathbb N_0^n } \frac 1 { \wp ! }
\partial_{ { y } }^\wp { f }_l \left( t, \varphi_0( t, c_0 )
\right) \tilde { { \delta } }^{ \, \wp }(t).
                                                        \label{TSE:NL}
\end{equation}

Подставляя формальное равенство \eqref{UFSolution} в векторное
дифференциальное уравнение \eqref{DEV:main} и учитывая формулы
\eqref{UFSolution'}, \eqref{TSE:NL}, получим новое формальное равенство
\begin{multline*}
\frac{ d \varphi_0( t, c_0 ) }{ dt } +
\sum_{ s=1 }^\infty \sum_{ p=1 }^{ \varkappa_s } \left(
\frac{ d \nu_{ sp }(t) }{ dt } \varphi_{ sp }( t ) + \nu_{ sp }(t)
\frac{ d \varphi_{ sp }(t) }{ dt } \right) =
                                                                        \\
= A \left( \varphi_0( t, c_0 ) + \sum_{ s=1 }^\infty
\sum_{ p=1 }^{ \varkappa_s } \nu_{ sp }(t) \varphi_{ sp }( t ) \right)
+ { f }(t) + \sum_{ l=1 }^m \varepsilon_l(t)
\sum_{ \wp \in \mathbb N_0^n } \frac 1 { \wp ! }
\partial_{ { y } }^\wp { f }_l \left( t, \varphi_0( t, c_0 )
\right) \tilde { { \delta } }^{ \, \wp }(t).
\end{multline*}
Перепишем это равенство в следующем виде
\begin{multline}
\frac{ d \varphi_0( t, c_0 ) }{ dt } - A \varphi_0( t, c_0 ) +
\sum_{ s=1 }^\infty \sum_{ p=1 }^{ \varkappa_s }
\nu_{ sp }(t) \left( \frac{ d \varphi_{ sp }(t) }{ dt } -
A \varphi_{ sp }( t ) \right) =
                                                                        \\
= { f }(t) - \sum_{ s=1 }^\infty \sum_{ p=1 }^{ \varkappa_s }
\frac{ d \nu_{ sp }(t) }{ dt } \varphi_{ sp }( t ) +
\sum_{ l=1 }^m \varepsilon_l(t)
\sum_{ \wp \in \mathbb N_0^n } \frac 1 { \wp ! }
\partial_{ { y } }^\wp { f }_l \left( t, \varphi_0( t, c_0 )
\right) \tilde { { \delta } }^{ \, \wp }(t).
                                                        \label{FEquality}
\end{multline}

Рассмотрим двойные суммы в правой части этого формального равенства.
В силу формулы \eqref{E:nu_sp(t)}, получим
\begin{multline*}
\frac{ d \nu_{ sp }(t) }{ dt } = \sum_{ k=1 }^m
\left( \prod^m_{ \substack{ l=1 \\ l \ne k } }
\prod_{ r=0 }^{ [ \rho_s - \varrho_l ] }
\left( \frac{ d^r \varepsilon_l(t) }{ d t^r } \right)^{ \beta_{ splr } }
\right) \sum_{ q=0 }^{ [ \rho_s - \varrho_k ] }
\left( \prod^{ [ \rho_s - \varrho_k ] }_{ \substack{ r=0 \\
r \ne q } } \left( \frac{ d^r \varepsilon_k(t) }{ d t^r } \right)^{
\beta_{ spkr } } \right) \beta_{ spkq } \cdot
                                                                        \\
\cdot \left( \frac{ d^q \varepsilon_k(t) }{ d t^q }
\right)^{ \beta_{ spkq } - 1 }
\frac{ d^{ q+1 } \varepsilon_k(t) }{ d t^{ q+1 } } =
\sum_{ k=1 }^m \sum_{ q=0 }^{ [ \rho_s - \varrho_k ] }
\beta_{ spkq } \left( \prod^m_{ \substack{ l=1 \\
l \ne k } } \prod_{ r=0 }^{ [ \rho_s - \varrho_l ] }
\left( \frac{ d^r \varepsilon_l(t) }{ d t^r } \right)^{ \beta_{ splr } }
\right) \cdot
                                                                        \\
\cdot \left( \prod^{ [ \rho_s - \varrho_k ] }_{
\substack{ r=0 \\ r \ne q } } \left(
\frac{ d^r \varepsilon_k(t) }{ d t^r } \right)^{ \beta_{ spkr } } \right)
\left( \frac{ d^q \varepsilon_k(t) }{ d t^q } \right)^{ \beta_{ spkq } -
1 } \frac{ d^{ q+1 } \varepsilon_k(t) }{ d t^{ q+1 } }
\ \ \left( s \in \mathbb N, \ p = \overline{ 1, \varkappa_s } \right).
\end{multline*}
Таким образом представили функции $ \dfrac{ d \nu_{ sp }(t) }{ dt } $
($ s \in \mathbb N $, $ p = \overline{ 1, \varkappa_s } $) в виде линейной
комбинации некоторых известных функций вида \eqref{E:nu_sp(t)} с
некоторыми известными действительными коэффициентами.
Причём ранг этих функций будет
\begin{multline}
\sum^m_{ \substack{ l=1 \\ l \ne k } }
\sum_{ r=0 }^{ [ \rho_s - \varrho_l ] } ( \varrho_l + r )
\beta_{ splr } + \sum^{ [ \rho_s - \varrho_k ] }_{
\substack{ r=0 \\ r \ne q } } ( \varrho_k + r )
\beta_{ spkr } + ( \varrho_k + q ) ( \beta_{ spkq } - 1 ) +                 \\
+ \varrho_k + q + 1 = 1 + \rho_s
\ \ \left( s \in \mathbb N, \ p = \overline{ 1, \varkappa_s } \right).
                                                        \label{Rank_nu'}
\end{multline}

Учитывая аксиомы класса $ \mathcal K_1(A) $ можем заключить, что произведение
двух любых компонент вектор--функции $ \tilde { { \delta } }(t) $ можно
записать в виде ряда, такого же, как и у этих множителей.
Поэтому правая часть равенства \eqref{TSE:NL} обладает формой, такой же, как
и у вектор--функции $ \tilde \delta(t) $.
Так как в этой части присутствуют функции $ \varepsilon_l(t) $
($ l = \overline{ 1, m} $), в качестве множителей, то в слагаемых ряда типа
$ \tilde \delta(t) $, представляющего правую часть равенства \eqref{TSE:NL},
при функциях $ \nu_{ sp }(t) $ ($ s \in \mathbb N $,
$ p = \overline{ 1, \varkappa_s } $) не будет вектор--функций
$ { \varphi }_{ kp }(t) $ ($ k \ge s $ ($ k \in \mathbb N $),
$ p = \overline{ 1, \varkappa_k } $).

Приравнивая в равенстве \eqref{FEquality} слева и справа слагаемые не
содержащие функций $ \nu_{ sp }(t) $ ($ s \in \mathbb N $,
$ p = \overline{ 1, \varkappa_s } $), получим векторное равенство.
В силу леммы \ref{L:truncated_main} вектор--функция
$ \varphi_0( t, c_0 ) \in \mathcal K_1^{ n \times 1 }(A) \cap C^1( I ) $
обращает его в тождество на промежутке $ I $.

Далее, приравнивая в равенстве \eqref{FEquality} с обеих сторон коэффициенты
при одинаковых функциях $ \nu_{ 1p }(t) $
($ p = \overline{ 1, \varkappa_1 } $),
получим совокупность векторных дифференциальных уравнений
\begin{equation}
\left[
\begin{array}{l}
\dfrac{ d \varphi_{ 1p }(t) }{ dt } = A { \varphi }_{ 1p }(t) +
{ f }_p ( t, \varphi_0( t, c_0 ) ),                                  \\
p = \overline{ 1, \varkappa_1 }.
\end{array}
\right.
                                                \label{PDEV:first_terms}
\end{equation}
В силу леммы \ref{L:truncated_main} у этих векторных дифференциальных
уравнений на промежутке $ I $ существует, по крайней мере, по одному частному
решению из класса $ \mathcal K_1^{ n \times 1 }(A) \cap C^1( I ) $ ---
вектор--функции $ { \varphi }_{ 1p }( t, c_{ 1p } ) $
($ p = \overline{ 1, \varkappa_1 } $, параметры
$ c_{ 1p } \in \mathbb R^{ n_{ 1p } } $, $ 0 \le n_{ 1p } \le n $).

Пусть уже найдены вектор--функции
$ \varphi_{ sp }( t, c_{ sp } ) \in \mathcal K_1^{ n \times 1 }(A) \cap C^1( I ) $
($ s = \overline{ 1, q } $, $ p = \overline{ 1, \varkappa_{ s } } $,
$ q \in \mathbb N $ фиксировано, параметры
$ c_{ sp } \in \mathbb R^{ n_{ sp } } $, $ 0 \le n_{ sp } \le n $).
Приравнивая в равенстве \eqref{FEquality} слева и справа коэффициенты при
одинаковых функциях $ \nu_{ q+1 \: p }(t) $
($ p = \overline{ 1, \varkappa_{ q+1 } } $),
получим систему векторных дифференциальных уравнений относительно пока
неизвестных вектор--функций $ { \varphi }_{ q+1 \: p }(t) $
($ p = \overline{ 1, \varkappa_{ q+1 } } $).
В силу свойств \eqref{Rank_nu'} и замечания о правой части равенства
\eqref{TSE:NL} правые части этих векторных дифференциальных уравнений уже
известны из предыдущих шагов.
Поэтому, учитывая вид левой части равенства \eqref{FEquality}, можем
заключить, что эти векторные дифференциальные уравнения являются линейными,
а сама система распадается на совокупность из $ \varkappa_{ q+1 } $
независимых уравнений.
В силу условий теоремы и аксиом класса $ \mathcal K_1(A) $ свободные
члены этих векторных дифференциальных уравнений принадлежат
$ \mathcal K_1^{ n \times 1 }(A) $.
Следовательно, снова в силу леммы \ref{L:truncated_main} у этих векторных
дифференциальных уравнений существует на промежутке $ I $, по крайней мере,
по одному частному решению из класса
$ \mathcal K_1^{ n \times 1 }(A) \cap C^1( I ) $ --- вектор--функции
$ \varphi_{ q+1 \: p }( t, c_{ q+1 \: p } ) $
($ p = \overline{ 1, \varkappa_{ q+1 } } $, параметры
$ c_{ q+1 \: p } \in \mathbb R^{ n_{ q+1 \: p } } $,
$ 0 \le n_{ q+1 \: p } \le n $).

Таким образом, с помощью математической индукции доказали, что формальное равенство
\eqref{FEquality} сводится к бесконечной последовательности совокупностей
векторных дифференциальных уравнений
\begin{equation}
\begin{array}{l}
\left[
\begin{array}{l}
\dfrac{ d \varphi_{ sp }(t) }{ dt } = A \varphi_{ sp }( t ) +
g_{ sp } \left( t, \varphi_0( t, c_0 ),
{ \varphi }_{ 11 }( t, c_{ 11 } ), \dots,
\varphi_{ s-1 \: \varkappa_{ s-1 } }( t, c_{ s-1 \: \varkappa_{ s-1 } } )
\right),
\\ p = \overline{ 1, \varkappa_s }                                        \\
\end{array}
\right.                                                                 \\
(s \in \mathbb N),
\end{array}
                                                        \label{SSDE:Coef}
\end{equation}
где вектор--функции $ g_{ sp } \left( t, \varphi_0( t, c_0 ),
{ \varphi }_{ 11 }( t, c_{ 11 } ), \dots,
\varphi_{ s-1 \: \varkappa_{ s-1 } }( t, c_{ s-1 \: \varkappa_{ s-1 } } )
\right) \in \mathcal K_1^{ n \times 1 }(A) $
($ s \in \mathbb N $, $ p = \overline{ 1, \varkappa_s } $) известны из
предыдущих шагов.
Причём из этих уравнений можно найти вектор--функции
$ \varphi_{ sp }( t, c_{ sp } ) $ ($ s \in \mathbb N $,
$ p = \overline{ 1, \varkappa_s } $, параметры
$ c_{ sp } \in \mathbb R^{ n_{ sp } } $, $ 0 \le n_{ sp } \le n $)
с требуемыми свойствами.

Как известно, общее решение системы обыкновенных дифференциальных уравнений
$ n $--го порядка зависит от $ n $ произвольных параметров, следовательно в
формальном семействе решений \eqref{FSolution} могут фигурировать не более
чем $ n $ независимых скалярных постоянных.

\end{proof}

\begin{remark}

Аналогичное теореме \ref{T:FSolution} утверждение можно получить для
случая, когда $  y : I \to \mathbb C^{ n \times 1 } $,
матрица $ A \in \mathbb C^{ n \times n } $,
$ { f }(t) : I \to \mathbb C^{ n \times 1 } $,
функции $ \varepsilon_l(t) : I \to \mathbb C $,
вектор--функции
$ { f }_l( t, { y } ) : D \to \mathbb C^{ n \times 1 } $
($ l = \overline{ 1, m } $), множество $ D \subset I \times \mathbb C^n $.

\end{remark}

\begin{remark}
                                                    \label{R:first_terms}
Векторные коэффициенты первых членов формального частного решения (семейства
решений) \eqref{FSolution} основного векторного дифференциального уравнения
\eqref{DEV:main} можно найти из совокупности векторных дифференциальных
уравнений \eqref{PDEV:first_terms}.

\end{remark}

\section{Асимптотический характер формального частного решения (семейства
решений)}
                                                    \label{CH:AsymptChart}

\subsection{Случай не чисто мнимых собственных значений матрицы A}

В этой части статьи исследован асимптотический характер формального частного
решения (семейства решений) \eqref{FSolution} основного вещественного
векторного квазилинейного обыкновенного дифференциального уравнения
\eqref{DEV:main} в том случае, когда матрица $ A $ не имеет чисто мнимых
собственных значений (теорема \ref{T:NPI_AsymptChart_O}).

\subsubsection{Вспомогательная лемма о выделении малости}

В следующей лемме после замены неизвестной вектор--функции вида
\eqref{CV:r_j} в основном векторном дифференциальном уравнении
\eqref{DEV:main}, выделена малость в свободном члене полученного в результате
этой замены векторного уравнения.

\begin{lemma}

Пусть для векторного дифференциального уравнения \eqref{DEV:main} выполнены
условия \ref{C:main_t_first}. -- \ref{C:main_t_second}., заданы некоторые
ранги $ \varrho_l $ функций $ \varepsilon_l(t) $, для которых выполняется
определение \ref{D:ProdDerivRank}, числа $ k \in \mathbb N $,
$ p_{ kl } := \max \left\{ 0,
\ \left[ \rho_k - \varrho_l \right] + 1 \right\} $
и
$ q_{ kl } := \left[ \frac { \rho_k - \varrho_l } { \varrho_1 } \right] $
($ l = \overline{ 1, m } $),
а также справедливы следующие предположения:
                                                        \label{L:separatedFT}
\begin{enumerate}
\item
функции $ \varepsilon_l(t) : I \to \mathbb R $,
$ \varepsilon_l(t) \in C^{ p_{ kl } } ( I ) $,
$ \dfrac{ d^r \varepsilon_l(t) }{ dt^r } =
O \left( \left( \varepsilon_1(t) \right)^{
\frac { \varrho_l + r } { \varrho_1 } } \right) $,
$ \varepsilon_1(t) \ge 0 $, $ \varepsilon_1(t) = o(1) $
($ l = \overline{ 1, m } $, $ r = \overline{ 0, p_{ kl } } $,
$ t \to + \infty $);
\item
вектор--функции
$ f_l( t, { y } ) : D \to \mathbb R^{ n \times 1 } $,
$ { f }_l( t, { y } ) \in C_{ t, \ { y } }^{ 0,
\ \max \left\{ 2, \ q_{ kl } + 1 \right\} } ( D ) $
($ l = \overline{ 1, m } $);
\item
супремумы
$ \sup\limits_{ ( t,  y ) \in D } \left\|
\partial_{ { y } }^\wp { f }_l( t, { y } ) \right\| < + \infty $
$ \left( \wp \in \mathbb N_0^n, \ | \wp | = \overline{ 0, \ \max \left\{ 2,
\ q_{ kl } + 1 \right\} }, \ l = \overline{ 1, m } \right) $;
\item
вектор--функции
$ \partial_{ { y } }^\wp { f }_l( t, \varphi_0( t, c_0 ) ) \in
\mathcal K_1^{ n \times 1 }(A) $
$ \left( \wp \in \mathbb N_0^n, \ | \wp | = \overline{ 0, q_{ kl } },
\ \rho_k \ge \varrho_l \right) $.
                                                        \label{C:F_l^(r) inaca}
\end{enumerate}
Тогда для достаточно большого числа $ t_0 $ с помощью замены неизвестной
вектор--функции вида
\begin{equation}
y = s( t ) + r,
                                                        \label{CV:r_j}
\end{equation}
где $ s( t ) $ --- усеченная сумма формального частного решения
(семейства решений) \eqref{FSolution} векторного дифференциального уравнения
\eqref{DEV:main}:
\begin{equation}
s( t ) = \varphi_0( t, c_0 ) + \sum_{ s=1 }^k
\sum_{ p=1 }^{ \varkappa_s } \nu_{ sp }(t)
\varphi_{ sp }( t, c_{ sp } ),
                                                        \label{s_j}
\end{equation}
вводя новую неизвестную вектор--функцию
$ r : I \to \mathbb R^{ n \times 1 } $, векторное дифференциальное уравнение
\eqref{DEV:main} можно привести к эквивалентному дифференциальному уравнению
вида
\begin{gather}
\frac{ d r }{ dt } = \hat A(t) r + g(t) +
\sum_{ l=1 }^m \varepsilon_l(t) { \psi }_l( t, r )
\ \ \left( \left( t, r \right) \in D_1 \right),
                                                    \label{DEV:separatedFT}
                                                                        \\
\hat A(t) := A + \sum_{ l=1 }^m \varepsilon_l(t)
\frac { \partial f_l( t, s( t ) ) } { \partial y },
                                                                    \notag
                                                                        \\
C(I) \ni  g(t) = O \left( \left( \varepsilon_1(t) \right)^{
\frac { \rho_{ k+1 } } { \varrho_1 } } \right)
\ \ \left( t \to + \infty \right),
                                                    \label{Est:FreeTerm}
\end{gather}
причём вектор--функции
$ \psi_l( t, r ) \in C_{ t, \ r }^{ 0, \ 2 } ( D_1 ) $
$ \left( l = \overline{ 1, m } \right), $ справедливы тождества
\begin{equation}
\psi_l( t, 0 ) \equiv 0 \ \ \left( l = \overline{ 1, m }, \ t \in I \right),
                                                    \label{Id:psi_l(t,0)=0}
\end{equation}
нормы
\begin{equation}
\left\| \frac{ \partial \psi_l(t, r) }{ \partial r }
\right\| = O \left( || r || \right)
\ \ \left( l = \overline{ 1, m }, \ t \to + \infty, \ r \to 0 \right),
                                                \label{Est:JM_psi_l(t,r)}
\end{equation}
область
$ D_1 := \left\{ ( t, r ) \ \left| \ t \in I,
\ r \in \mathbb R^{ n \times 1 }, \ || r || \le \hat a_0 \right. \right\}, $
число $ \hat a_0 < a $ ($ \hat a_0 \in \mathbb R^+ $, разность
$ a - \hat a_0 $ можно сделать сколь угодно малой за счёт выбора числа
$ t_0 $ достаточно большим).

\end{lemma}

\begin{proof}
Поступая точно также как и в теореме \ref{T:FSolution} для нахождения
вектор--функций
$ \varphi_{ sp }( t, c_{ sp } ) \in
\mathcal K_1^{ n \times 1 }(A) \cap C^1( I ) $
($ s = \overline{ 1, k } $, $ p = \overline{ 1, \varkappa_s } $), в силу
условий леммы, найдем вектор--функцию $ s( t ) \in C^1(I) $.
(При этом свободные члены векторных дифференциальных уравнений
\eqref{SSDE:Coef} могут содержать только вектор--функции
$ \partial_{ { y } }^\wp { f }_l( t, \varphi_0( t, c_0 ) ) $,
$ \wp \in \mathbb N_0^n $, $ | \wp | = \overline{ 0, q_{ kl } } $,
$ \rho_k \ge \varrho_l $).

Сделаем в векторном дифференциальном уравнении \eqref{DEV:main} замену
неизвестной вектор--функции вида \eqref{CV:r_j}, в результате получим
дифференциальное уравнение вида
\begin{equation}
\frac{ d s( t ) }{ dt } + \frac{ dr }{ dt } =
A ( s( t ) + { r } ) + { f }(t) +
\sum_{ l=1 }^m \varepsilon_l(t) f_l \left( t, s( t ) + r \right).
                                                        \label{SDE:mainr_j}
\end{equation}
С помощью формулы Тейлора получим разложения вектор--функций
\begin{equation}
{ f }_l ( t, s( t ) + { r } ) =
{ f }_l ( t, s( t ) ) +
\frac { \partial { f }_l ( t, s( t ) ) }
{ \partial { y } } { r } + { \psi }_l ( t, { r } )
\ \ ( l = \overline{ 1, m } ),
                                                        \label{psi_l(t,r)}
\end{equation}
причём вектор--функции
$ \psi_l ( t, r ) \in C_{ t, \ r }^{ 0, \ 2 } ( D_1 ), $
$ \psi_l ( t, r ) = o( || r || ) $ при $ t \to + \infty, $ $ r \to 0 $
($ l = \overline{ 1, m } $).
Учитывая эти свойства, перепишем векторное дифференциальное уравнение
\eqref{SDE:mainr_j} в виде
\begin{gather}
\frac{ dr }{ dt } = A { r } +  g(t) +
\sum_{ l=1 }^m \varepsilon_l(t) \left(
\frac { \partial { f }_l ( t, s( t ) ) }
{ \partial { y } } { r }
+ { \psi }_l ( t, { r } ) \right),
                                                    \label{SDE:separatedsr}
                                                                        \\
g(t) := - \frac{ d s( t ) }{ dt } + A s( t ) +
{ f }(t) + \sum_{ l=1 }^m \varepsilon_l(t) { f }_l ( t, s( t ) ).
                                                        \label{E:FreeTerm}
\end{gather}

Оценим свободный член $  g(t) $ последнего векторного дифференциального
уравнения.
Для этого в равенстве \eqref{E:FreeTerm} заменим некоторые вектор--функции
$ { f }_l ( t, s( t ) ) $ их разложениями по формуле Тейлора и докажем, что в
правой части формулы \eqref{E:FreeTerm} после этого, за счёт
взаимоуничтожений, останутся только слагаемые содержащие функции вида
\eqref{E:nu_sp(t)}, ранг которых больше чем $ \rho_k $.
А затем, оценим полученное выражение с помощью ранга.
Итак, заменим в равенстве \eqref{E:FreeTerm} вектор--функцию
$ s( t ) $ по формуле \eqref{s_j}.
Для $ \varrho_l > \rho_k $ заменим в \eqref{E:FreeTerm} вектор--функции
$ { f }_l ( t, s( t ) ) $ их оценкой --- $ O (1) $
($ t \to + \infty $).
Чтобы рассмотреть случай, когда $ \varrho_l \le \rho_k $ нам потребуются
нижеследующие оценки функций.

В силу определения \ref{ProdDerivRank} и условий леммы имеют место оценки
\begin{equation}
\nu_{ sp }(t) = O \left( \left( \varepsilon_1(t)
\right)^{ \frac { \rho_s } { \varrho_1 } } \right)
\ \ \left( t \to + \infty, \ s \in \mathbb N,
\ p = \overline{ 1, \varkappa_s } \right),
                                                        \label{estimate_nu}
\end{equation}
где в представлениях \eqref{E:nu_sp(t)} функций $ \nu_{ sp }(t) $
($ s \in \mathbb N $, $ p = \overline{ 1, \varkappa_s } $), могут
присутствовать производные функций $ \varepsilon_l(t) $
($ l = \overline{ 1, m } $) порядка не выше чем $ p_{ kl } $
($ l = \overline{ 1, m } $).
Следовательно, вектор--функция
\begin{equation}
\delta( t ) := s( t ) - \varphi_0( t, c_0 ) =
\varepsilon_1(t) O( 1 ) \ \ ( t \to + \infty ).
\end{equation}
А значит, функции
\begin{equation}
\delta^\wp( t ) = \varepsilon_1^{ | \wp | }(t) O( 1 )
\ \ ( t \to + \infty, \ \wp \in \mathbb N_0^n, \ \wp \ne 0 ).
                                                \label{estimate_delta_K}
\end{equation}

Учитывая последние оценки, с помощью формулы Тейлора получим разложения
вектор--функций
\begin{equation*}
f_l( t, s( t ) ) =
\sum_{ \substack{ 0 \le | \wp | \le q_{ kl }
\\ ( \wp \in \mathbb N_0^n) } } \frac 1 { \wp ! }
\partial_{ { y } }^\wp { f }_l ( t, \varphi_0( t, c_0 ) ) \delta^\wp( t ) +
\left( \varepsilon_1(t) \right)^{ q_{ kl } + 1 } O( 1 )
\ \ ( t \to + \infty, \ \varrho_l \le \rho_k ).
\end{equation*}
Заменим вектор--функции
$ { f }_l ( t, s( t ) ) $ ($ \varrho_l \le \rho_k $) в равенстве
\eqref{E:FreeTerm} этими разложениями.
В итоге, учитывая оценки \eqref{estimate_nu}, формула \eqref{E:FreeTerm}
примет вид
\begin{multline}
 g(t) = - \frac{ d \varphi_0( t, c_0 ) }{ dt } + A \varphi_0( t, c_0 ) -
\sum_{ s=1 }^k \sum_{ p=1 }^{ \varkappa_s } \nu_{ sp }(t)
\left( \frac{ d \varphi_{ sp }( t, c_{ sp } ) }{ dt } -
A \varphi_{ sp }( t, c_{ sp } ) \right) + { f }(t) -
                                                                        \\
- \sum_{ \substack{ \rho_s + 1 \le \rho_k  \\ ( s \in \mathbb N ) } }
\sum_{ p=1 }^{ \varkappa_s } \frac{ d \nu_{ sp }(t) }{ dt }
\varphi_{ sp }( t, c_{ sp } ) +
\sum_{ \varrho_l \le \rho_k } \varepsilon_l(t)
\sum_{ \substack{ 0 \le | \wp | \le q_{ kl }
\\ ( \wp \in \mathbb N_0^n ) } } \frac 1 { \wp ! }
\partial_{ { y } }^\wp { f }_l ( t, \varphi_0( t, c_0 ) ) \delta^\wp( t ) +
                                                                        \\
+ \sum_{ \substack{ \rho_s + 1 > \rho_k  \\
( s \in \mathbb N, \ s \le k ) } }
\sum_{ p=1 }^{ \varkappa_s } \frac{ d \nu_{ sp }(t) }{ dt } O ( 1 ) +
\sum_{ \varrho_l \le \rho_k } \varepsilon_l(t)
\left( \varepsilon_1(t) \right)^{ q_{ kl } + 1 } O ( 1 ) +
\sum_{ \varrho_l > \rho_k } \varepsilon_l(t) O ( 1 )
\ \ ( t \to + \infty ).
                                                    \label{E:FreeTermExpan}
\end{multline}

Покажем, что разложения вектор--функций
$ f_l( t, s( t ) ) $ ($ \varrho_l \le \rho_k $)
содержат достаточное количество слагаемых.
Для этого оценим снизу
\begin{multline*}
\operatorname{ Rank } \left( \varepsilon_l(t)
\left( \varepsilon_1(t) \right)^{ q_{ kl } + 1 } \right) =
\varrho_l + \varrho_1 \left( q_{ kl } + 1 \right) =
\varrho_l + \varrho_1 \left( \frac { \rho_k - \varrho_l }
{ \varrho_1 } - \left\{ \frac { \rho_k - \varrho_l }
{ \varrho_1 } \right\} + 1 \right) =
                                                                        \\
= \rho_k + \varrho_1 \left( 1 - \left\{ \frac { \rho_k -
\varrho_l } { \varrho_1 } \right\} \right) > \rho_k
\ \ (\varrho_l \le \rho_k).
\end{multline*}

Таким образом, ранги коэффициентов при $ O $ у всех слагаемых из равенства
\eqref{E:FreeTermExpan}, содержащих этот символ, строго больше чем
$ \rho_k $.
Так как ранг может принимать только дискретные значения из множества $ R $,
то ранги этих коэффициентов будут не меньше чем $ \rho_{ k+1 } $.

Далее, если в формальном равенстве \eqref{FEquality} все члены перенести в
правую часть, неизвестные вектор--функции $ \varphi_{ sp }( t ) $ заменить на
уже найденные $ \varphi_{ sp }( t, c_{ sp } ) $ ($ s = \overline{ 1, k } $,
$ p = \overline{ 1, \varkappa_s } $) и исключить из неё все слагаемые
содержащие функции вида \eqref{E:nu_sp(t)}, ранги которых больше чем
$ \rho_k $, то полученное выражение будет совпадать с правой частью формулы
\eqref{E:FreeTermExpan}, если из неё тоже исключить такие же слагаемые.
Поэтому в силу определения вектор--функций $ \varphi_0( t, c_0 ) $ и
$ \varphi_{ sp }( t, c_{ sp } ) $, $ \nu_{ sp }(t) $ и условий леммы свойство
\eqref{E:FreeTermExpan} можно переписать в виде конечной суммы
$$
 g(t) = \sum_{
\substack{ s \ge k + 1  \\ ( s \in \mathbb N ) } }
\sum_{ p=1 }^{ \varkappa_s } \nu_{ sp }(t) O(1)
\ \ ( t \to + \infty ),
$$
причём в представлениях функций $ \nu_{ sp }(t) $ могут фигурировать
производные функций $ \varepsilon_l(t) $ ($ l = \overline{ 1, m } $)
порядка не больше чем $ p_{ kl } $ ($ l = \overline{ 1, m } $).
На основании оценок \eqref{estimate_nu} имеем право переписать это равенство
в виде \eqref{Est:FreeTerm}.

В итоге, векторное дифференциальное уравнение \eqref{SDE:separatedsr}
привели к виду \eqref{DEV:separatedFT}.

В силу разложений \eqref{psi_l(t,r)}, очевидно, что выполняются тождества
\eqref{Id:psi_l(t,0)=0}.

Оценим сверху нормы матриц Якоби вектор--функций $ \psi_l ( t, r ) $ при
$ t \to + \infty $, $ r \to 0 $ ($ l = \overline{ 1, m } $).
Дифференцируя по $ r $ частным образом разложения \eqref{psi_l(t,r)} и
перенося слагаемые, получим выражения матриц Якоби вектор--функций
$ \psi_l( t, r ) $ через матрицы Якоби вектор--функций $ f_l( t, y ) $
\begin{equation*}
\frac{ \partial \psi_l( t, r ) }{ \partial r } =
\frac{ \partial f_l( t, s( t ) + r ) }{ \partial y } -
\frac{ \partial f_l( t, s( t ) ) }{ \partial y }
\ \ \left( l = \overline{ 1, m } \right).
\end{equation*}
Оценим сверху нормы этих разностей матриц Якоби.
Для этого применим к каждой строке правой части последних матричных равенств
аналог формулы Лагранжа конечных приращений скалярных функций, в результате
получим неравенства
\begin{multline*}
\left\| \nabla_{ y } \left( f_l( t, s( t ) + r ) \right)_j
- \nabla_{ y } \left( f_l( t, s( t ) ) \right)_j \right\| \le
\sup_{ ( t, y ) \in D } \left\|
\frac{ \partial \nabla_{ y } \left( f_l( t, y ) \right)_j }
{ \partial y } \right\| O \left( || r || \right)
\\ \left( l = \overline{ 1, m }, \ j = \overline{ 1, n },
\ t \to + \infty, \ r \to 0 \right).
\end{multline*}
Следовательно, нормы
\begin{multline*}
\left\| \frac{ \partial f_l( t, s( t ) + r ) } { \partial y } -
\frac{ \partial f_l( t, s( t ) ) }{ \partial y } \right\| \le
\max_{ j = \overline{ 1, n } } \sup_{ ( t, y ) \in D } \left\|
\frac{ \partial \nabla_{ y } \left( f_l( t, y ) \right)_j }
{ \partial y } \right\| O \left( || r || \right)
                                                                        \\
\left( l = \overline{ 1, m }, \ t \to + \infty, \ r \to 0 \right).
\end{multline*}
В силу условий леммы супремумы в последних неравенствах будут ограничены,
поэтому имеют силу оценки \eqref{Est:JM_psi_l(t,r)}.

\end{proof}

\begin{remark}

Далее будем считать, что свободный член $ g(t) $ в векторном дифференциальном
уравнении \eqref{DEV:separatedFT} отличен от тождественного нуля на
промежутке $ I $.
Так как в противном случае вектор--функция $ r \equiv 0 $ --- точное частное
решение этого векторного дифференциального уравнения на промежутке $ I $.
А вектор--функция $ y \equiv s( t ) $ --- точное частное решение (семейство
решений) основного векторного дифференциального уравнения \eqref{DEV:main} на
промежутке $ I $.

\end{remark}

\subsubsection{Теорема об асимптотическом характере}

В ходе доказательства теоремы \ref{T:NPI_AsymptChart_O} будет использовано
следующее понятие.

\begin{definition}

Пусть $ f(t) $ --- вещественная или комплексная функция, непрерывная при всех
конечных значениях $ t \ge 0 $.
Будем говорить, что $ f(t) $ --- функция со слабой вариацией (смотрите
\cite{Bib:Persidskiy}), если для любого наперед заданного сколь угодно малого
числа $ \epsilon > 0 $ и для любого наперед заданного сколь угодно большого
числа $ \alpha $ существует такое число
$ \beta = \beta( \epsilon, \alpha ), $ что при любых значениях
$ t_1 \ge \beta $ и $ t_2 \ge \beta $ имеет место неравенство
$ | f( t_1 ) - f( t_2 ) | < \epsilon $, если только
$ | t_1 - t_2 | < \alpha. $

\end{definition}

\begin{example}
                                                \label{Ex:weak_variation}
Если при $ t \to + \infty $ функция $ f(t) $ стремится к некоторому конечному
пределу, то $ f(t) $ есть функция со слабой вариацией
(смотрите \cite{Bib:Persidskiy}).

\end{example}

В следующей теореме исследован асимптотический характер формального частного
решения (семейства решений) \eqref{FSolution} основного вещественного
векторного квазилинейного обыкновенного дифференциального уравнения
\eqref{DEV:main} в том случае, когда матрица $ A $ не имеет чисто мнимых
собственных значений.

\begin{theorem}
                                                \label{T:NPI_AsymptChart_O}
Пусть для векторного дифференциального уравнения \eqref{DEV:main} выполнены
все условия леммы \ref{L:separatedFT}, а также справедливы следующие
предположения:
\begin{enumerate}
\item
$ \Re \, \lambda_j(A) \ne 0 $ $ \left( j = \overline{ 1, n } \right) $;
\item
функция $ \varepsilon_1(t) > 0 $ $ ( t \in I ) $.
\end{enumerate}
Тогда для достаточно большого числа $ t_0 $ у векторного дифференциального
уравнения \eqref{DEV:main} на промежутке $ I $ существует хотя бы одно
частное решение вида \eqref{CV:r_j}, где вектор--функции
$ \varphi_0( t, c_0 ), \varphi_{ sp }( t, c_{ sp } ) \in
\mathcal K_1^{ n \times 1 }(A) \cap C^1( I ) $
($ s = \overline{ 1, k } $, $ p = \overline{ 1, \varkappa_s } $), причём
погрешность
\begin{equation*}
r = O \left( \left( \varepsilon_1(t) \right)^{
\frac { \rho_{ k+1 } } { \varrho_1 } } \right)
\ \ ( t \to + \infty ).
\end{equation*}
Более того, на промежутке $ [ t_1, + \infty ) $ погрешность $ r $ зависит от
стольких произвольных скалярных постоянных, сколько имеется индексов $ j $
таких, что справедливы условия
\begin{equation}
j \in \left\{ 1, 2, 3, \dots, n \right\}, \ \ \Re \, \lambda_j(A) < 0,
                                            \label{C:NPI_No_of_parameters_O}
\end{equation}
$ t_1 \in \mathbb R $ --- достаточно большое число, определяется этими
постоянными ($ t_1 \ge t_0 $).
(Семейство решений вида \eqref{CV:r_j} может зависеть не более чем от $ n $
произвольных скалярных параметров.)

\end{theorem}

\begin{proof}
В силу леммы \ref{L:separatedFT} для достаточно большого числа $ t_0 $ с
помощью замены неизвестной вектор--функции вида \eqref{CV:r_j}, вводя новую
неизвестную вектор--функцию $ r : I \to \mathbb R^{ n \times 1 } $, векторное
дифференциальное уравнение \eqref{DEV:main} можно привести к эквивалентному
дифференциальному уравнению \eqref{DEV:separatedFT}.

Добьемся того, чтобы отношения элементов находящихся по обе стороны от
главной диагонали матрицы линейной однородной части этого векторного
дифференциального уравнения к диагональным были сколь угодно малы по
абсолютной величине.
Для этого воспользуемся теоремой 9 (страница 26) из статьи \cite{Bib:Persidskiy}.
Так как для достаточно большого числа $ t_0 $ матрица $ \hat A( t ) $ состоит
из ограниченных функций со слабой вариацией (смотрите пример
\ref{Ex:weak_variation}), то линейное однородное векторное дифференциальное
уравнение соответствующее уравнению \eqref{DEV:separatedFT} удовлетворяет
всем требованиям этой теоремы.
В её силу существует такая матрица $ K(t) = K( t, d_0 ) $ ($ d_0 > 0 $ ---
наперед заданное число), $ K(t) : I \to \mathbb R^{ n \times n } $,
$ K(t) \in C^1(I) $,
$$
\sup\limits_{ t \in I } || K(t) || < + \infty,
\ \sup\limits_{ t \in I } \left\| K^{ -1 }(t) \right\| < + \infty,
\ \sup\limits_{ t \in I } \left\| \dfrac{ d K(t) }{ dt } \right\| < + \infty,
$$
что в результате линейной замены неизвестной вектор--функции вида
$ r = K(t) x $ в векторном дифференциальном уравнении
\eqref{DEV:separatedFT} получим уравнение вида
\begin{equation}
\dfrac{ dx }{ dt } = ( U(t) + H(t) ) x + K^{ -1 }(t) g(t) +
K^{ -1 }(t) \sum_{ l=1 }^m \varepsilon_l(t) \psi_l \left( t, K(t) x \right)
\ \ \left( \left( t, x \right) \in \hat D_2 \right),
                                                \label{DEV:After_Persidsky}
\end{equation}
где матрица $ U(t) := \operatorname{ diag } ( u_1(t), \dots, u_n(t) ) $,
функции $ u_j(t) $ $ \left( j = \overline{ 1, n } \right) $ являются
вещественными частями корней уравнения
$ \det \left( \hat A(t) - \lambda E \right) = 0 $ относительно $ \lambda $,
для достаточно большого числа $ t_0 $ супремум
$ \sup\limits_{ t \in I } || H(t) || < d_0 $, область
$$
\hat D_2 := \left\{ ( t, x ) \ \left| \ t \in I,
\ x \in \mathbb R^{ n \times 1 }, \ || x || \le \hat a_0 \left( n
\sup\limits_{ t \in I } || K(t) || \right)^{ -1 } \right. \right\}.
$$

Так как собственные значения квадратной комплексной матрицы непрерывно зависят
от её элементов и матрица $ \hat A(t) \in C(I) $,
$ \lim\limits_{ t \to + \infty } \hat A(t) = A $, то
$ u_j(t) = \Re \, \lambda_j(A) + o(1) $
$ \left( j = \overline{ 1, n }, \ t \to + \infty \right) $.

Выберем малость в свободном члене векторного дифференциального уравнения
\eqref{DEV:After_Persidsky}.
Для этого сделаем в нём замену неизвестной вектор--функции вида
$$
x = \left( \varepsilon_1(t) \right)^{ \frac { \rho_{ k+1 } } { \varrho_1 } }
z.
$$
В результате получим векторное дифференциальное уравнение вида
\begin{gather}
\dfrac{ dz }{ dt } = \hat P(t) z + \hat q(t) + \hat x( t, z )
\ \ \left( \left( t, z \right) \in \hat D_3 \right),
                                                \label{DEV:NPI_solution->0}
                                                                        \\
\hat P(t) := U(t) - \frac{ \rho_{ k+1 } }{ \varrho_1 } \left(
\varepsilon_1(t) \right)^{ -1 } \frac{ d \varepsilon_1(t) }{ dt } E + H(t),
                                                        \label{NPI_P(t)}
                                                                        \\
\hat q(t) := \left( \varepsilon_1(t) \right)^{ -
\frac{ \rho_{ k+1 } } { \varrho_1 } } K^{ -1 }(t) g(t),
                                                        \label{NPI_q(t)}
                                                                        \\
\hat x( t, z ) := \left( \varepsilon_1(t) \right)^{ -
\frac{ \rho_{ k+1 } }{ \varrho_1 } } K^{ -1 }(t) \sum_{ l=1 }^m
\varepsilon_l(t) \psi_l \left( t, \left( \varepsilon_1(t) \right)^{
\frac{ \rho_{ k+1 } } { \varrho_1 } } K(t) z \right),
                                                        \label{NPI_x(t,z)}
\end{gather}
где область
$ \hat D_3 := \left\{ ( t, z ) \ \left| \ t \in I,
\ z \in \mathbb R^{ n \times 1 }, \ || z || \le \hat a \right. \right\}, $
число
\[
\hat a := \hat a_0 \left( \sup\limits_{ t \in I } \varepsilon_1(t) \right)^{
- \frac { \rho_{ k+1 } } { \varrho_1 } }
\left( n \sup\limits_{ t \in I } || K(t) || \right)^{ -1 }.
\]

Проверим, выполнено ли для вектор-функции $ \hat x( t, z ) $ в области
$ \hat D_3 $ по переменной $ z $ условие Липшица.
Для оценки сверху нормы $ || \hat x( t, z_1 ) - \hat x( t, z_2 ) || $
$ \left( \forall \ ( t, z_1 ), ( t, z_2 ) \in \hat D_3 \right) $
воспользуемся аналогом формулы Лагранжа конечных приращений скалярных
функций.
Таким образом, достаточно оценить сверху норму матрицы Якоби вектор--функции
$ \hat x( t, z ) $ по переменной $ z $ в области $ \hat D_3 $.
Найдем эту матрицу, используя формулу \eqref{NPI_x(t,z)}:
$$
\frac{ \partial \hat x( t, z ) }{ \partial z } = K^{ -1 }(t)
\sum_{ l=1 }^m \varepsilon_l(t) \frac{ \partial \psi_l \left( t,
\left( \varepsilon_1(t) \right)^{ \frac { \rho_{ k+1 } } { \varrho_1 } } K(t)
z \right) } { \partial r } K(t).
$$
Следовательно,
$$
\left\| \frac{ \partial \hat x( t, z ) }{ \partial z } \right\| =
O( \varepsilon_1(t) ) \sum_{ l=1 }^m \left\| \frac{ \partial \psi_l \left( t,
\left( \varepsilon_1(t) \right)^{ \frac{ \rho_{ k+1 } }{ \varrho_1 } } K(t) z
\right) }{ \partial r } \right\| \ \ \left( t \to + \infty \right).
$$
В результате подстановки свойств \eqref{Est:JM_psi_l(t,r)} в это неравенство,
получим искомую оценку матрицы Якоби
$$
\left\| \frac{ \partial \hat x( t, z ) }{ \partial z } \right\| =
O \left( \left( \varepsilon_1(t)
\right)^{ 1 + \frac{ \rho_{ k+1 } }{ \varrho_1 } } \right)
\ \ \left( t \to + \infty, \ || z || \le \hat a \right).
$$
В итоге получили требуемое условие Липшица для вектор-функции
$ \hat x( t, z ) $ в области $ \hat D_3 $ по переменной $ z $:
\begin{multline}
|| \hat x( t, z_1 ) - \hat x( t, z_2 ) || \le
\hat x(t) || z_1 - z_2 ||,
\ \ 0 < \hat x(t) := const \left( \varepsilon_1(t) \right)^
{ 1 + \frac { \rho_{ k+1 } } { \varrho_1 } }
                                                                        \\
\left( \forall \ ( t, z_1 ), ( t, z_2 ) \in \hat D_3,
\ t \to + \infty \right).
                                                        \label{IE:NPI_Lip}
\end{multline}

В силу тождеств \eqref{Id:psi_l(t,0)=0} и равенства \eqref{NPI_x(t,z)},
очевидно, что вектор--функция $ \hat x( t, 0 ) \equiv 0 $ $ ( t \in I ) $.
Поэтому, учитывая условие Липшица \eqref{IE:NPI_Lip}, заметим, что для
нелинейности $ \hat x(t, z) $ в области $ \hat D_3 $ будет иметь место оценка
\begin{equation}
|| \hat x(t, z) || \le \hat x(t) || z || \ \ \left( t \to + \infty \right).
                                                \label{IE:NPI_estimateNL}
\end{equation}

Запишем векторное дифференциальное уравнение \eqref{DEV:NPI_solution->0} в
скалярной форме
\begin{gather}
\dfrac{ d z_j }{ dt } = \hat q_j(t) + \sum_{ l=1 }^n \hat p_{ jl }(t) z_l +
\hat x_j(t, z_1, \dots, z_n)
\ \ \left( j = \overline{ 1, n } \right),
                                                \label{SDE:NPI_solution->0}
                                                                        \\
z_j := ( z )_j, \ \hat q_j(t) := ( \hat q(t) )_j,
\ \hat p_{ jl }(t) := ( \hat P(t) )_{ jl },
\ \hat x_j(t, z_1, \dots, z_n) := ( \hat x( t, z ) )_j,
                                                                \notag
\end{gather}
где $ l = \overline{ 1, n }. $
Применим для нахождения ограниченных при $ t \ge t_0 $ решений системы
дифференциальных уравнений \eqref{SDE:NPI_solution->0} специальный метод
последовательных приближений, аналогичный тому, который использован в ходе
доказательства теоремы 1.1 из \S \ 1 главы IV кандидатской диссертации
\cite{Bib:Kostin_PHDT} (страница 67).
Пусть $ z_{ j \, s-1 }(t) $ ($ j = \overline{ 1, n } $) означает $ ( s-1 ) $-е
приближение, а $ z_{ js }(t) $ --- $ s $-е ($ j = \overline{ 1, n } $).
Положим
\[
z_{ 10 }(t) := \dots := z_{ n0 }(t) := 0.
\]
Определим $ s $-е приближение из системы
\begin{equation}
\dfrac{ d z_{ js } }{ dt } = \hat q_j(t) + \sum_{ l=1 }^{ j-1 }
\hat p_{ jl }(t) z_{ l \, s-1 } + \hat x_j(t, z_{ 1 \, s-1 }, \dots,
z_{ n \, s-1 } ) + \sum_{ l=j }^n \hat p_{ jl }(t) z_{ ls }
\ \ \left( j = \overline{ 1, n } \right),
                                        \label{SDE:NPI_Sth-approximation}
\end{equation}
выбирая начальные значения для функций $ z_{ js }(t) $
($ j = \overline{ 1, n } $) так, чтобы эти функции выражались через
$ z_{ j \, s-1 }(t) $ ($ j = \overline{ 1, n } $) формулами вида
\begin{multline}
z_{ js }(t) = \int\limits_{ \hat a_j }^t \hat q_j( \tau )
\exp \int\limits_\tau^t \hat p_{ jj }(t) \, dt \, d \tau +
\sum_{ l=1 }^{ j-1 } \int\limits_{ \hat a_j }^t \hat p_{ jl }( \tau )
z_{ l \, s-1 }( \tau ) \exp \int\limits_\tau^t \hat p_{ jj }(t)
\, dt \, d \tau +
                                                                        \\
+ z_j(t_0) \exp \int\limits_{ t_0 }^t \hat p_{ jj }( \tau ) \, d \tau +
\int\limits_{ \hat a_j }^t \hat x_j( \tau, z_{ 1 \, s-1 }( \tau ), \dots,
z_{ n \, s-1 }( \tau ) ) \exp \int\limits_\tau^t \hat p_{ jj }(t)
\, dt \, d \tau +
                                                                        \\
+ \sum_{ l = j+1 }^n \int\limits_{ \hat a_j }^t \hat p_{ jl }( \tau )
z_{ l s }( \tau ) \exp \int\limits_\tau^t \hat p_{ jj }(t) \, dt \, d \tau
\ \ \left( j = \overline{ 1, n } \right),
                                                    \label{NPI_solution->0}
\end{multline}
где каждый предел интегрирования $ \hat a_j $ равен либо $ t_0 $, либо
$ + \infty $; начальные значения $ z_j(t_0) := 0 $ для тех индексов
$ j $, для которых не выполняются условия \eqref{C:NPI_No_of_parameters_O}.
Нетрудно проверить простым дифференцированием, что при любом выборе указанных
пределов интегрирования, равенства \eqref{NPI_solution->0} будут давать нам
некоторое частное решение системы дифференциальных уравнений
\eqref{SDE:NPI_Sth-approximation} (если только величины, входящие в
\eqref{NPI_solution->0} не теряют смысла).
Пределы интегрирования $ \hat a_j $ будем выбирать так:
\[
\hat a_j :=
\begin{cases}
t_0, & \text{если } \Re \, \lambda_j(A) < 0,
                                                                        \\
+ \infty, & \text{если } \Re \, \lambda_j(A) > 0
\end{cases}
\ \ \left( j = \overline{ 1, n } \right).
\]

Постараемся сделать так, чтобы все последовательные приближения были
ограничены по модулю одним и тем же числом.
С этой целью предположим, что
\begin{equation}
| z_{ j \, s-1 }(t) | \le \hat \epsilon_0
\ \ ( j = \overline{ 1, n } ),
                                                \label{IE:NPI_constraint}
\end{equation}
где $ \hat \epsilon_0 $ --- некоторая константа,
$ 0 < \hat \epsilon_0 \le \hat a $ и потребуем, чтобы такие же неравенства
имели место и для $ s $-го приближения.
Принимая во внимание свойство \eqref{IE:NPI_estimateNL}, нетрудно заметить,
что при выполнении неравенств \eqref{IE:NPI_constraint}, функции
$ z_{ js }(t) $ ($ j = \overline{ 1, n } $) будут мажорироваться по модулю
функциями $ \hat \xi_j(t, t_0, \hat \epsilon_0) $
($ j = \overline{ 1, n } $), которые последовательно определяются из равенств
\begin{multline}
\hat \xi_j(t, t_0, \hat \epsilon_0) := \hat q_j^*(t) +
\hat \epsilon_0 \sum_{ l=1 }^{ j-1 } \hat p_{ jl }^*(t) + \hat z_j^*(t) +
\hat \epsilon_0 \hat x_j^*(t) +
                                                                        \\
+ \hat b_j \sum_{ l = j+1 }^n \int\limits_{ \hat a_j }^t \left|
\hat p_{ jl }( \tau ) \right| \hat \xi_l( \tau, t_0, \hat \epsilon_0)
\exp \int\limits_\tau^t \hat p_{ jj }(t) \, dt \, d \tau
\ \ ( j = \overline{ 1, n } ),
                                            \label{NPI_xi_j(t t_0 epsilon_0)}
\end{multline}
где функции
\begin{gather*}
\hat q_j^*(t) := \hat b_j \int\limits_{ \hat a_j }^t
\left| \hat q_j( \tau ) \right| \exp \int\limits_\tau^t \hat p_{ jj }(t)
\, dt \, d \tau,
\ \ \hat p_{ jl }^*(t) := \hat b_j \int\limits_{ \hat a_j }^t \left|
\hat p_{ jl }( \tau ) \right| \exp \int\limits_\tau^t
\hat p_{ jj }(t) \, dt \, d \tau,
                                                                        \\
\hat z_j^*(t) := | z_j(t_0) | \exp \int\limits_{ t_0 }^t
\hat p_{ jj }( \tau ) \, d \tau,
\ \ \hat x_j^*(t) := \hat b_j \int\limits_{ \hat a_j }^t \hat x( \tau )
\exp \int\limits_\tau^t \hat p_{ jj }(t) \, dt \, d \tau
\end{gather*}
$ \left( j, l = \overline{ 1, n }; \ j \ne l \right), $
числа
\[
\hat b_j :=
\begin{cases}
\phantom{ - } 1,      & \text{если } \Re \, \lambda_j(A) < 0,
                                                                        \\
-1,     & \text{если } \Re \, \lambda_j(A) > 0
\end{cases}
\ \ \left( j = \overline{ 1, n } \right).
\]

Для оценки сверху функций $ \hat \xi_j(t, t_0, \hat \epsilon_0) $
$ \left( j = \overline{ 1, n } \right) $ поступим аналогично тому, как это
сделано в доказательстве теоремы 2.1 из \S \ 2 главы IV кандидатской
диссертации \cite{Bib:Kostin_PHDT} (страница 69).
Докажем, что за счёт выбора достаточно малого числа $ d_0 $ и достаточно
большого $ t_0 $ можем сделать супремумы
$ \sup\limits_{ t \in I } \hat q_j^*(t) $,
$ \sup\limits_{ t \in I } \hat p_{ jl }^*(t) $,
$ \sup\limits_{ t \in I } \hat x_j^*(t) $
$ \left( j, l = \overline{ 1, n }; \ j \ne l \right) $ сколь угодно малыми.
Учитывая введённые ранее обозначения \eqref{NPI_P(t)}, заметим, что для
достаточно малого числа $ d_0 $ и достаточно большого числа $ t_0 $ функции
$ \hat p_{ jj }(t) $ сохраняют знак в строгом смысле на промежутке $ I $,
интегралы
\[
\int\limits_I \hat p_{ jj }(t) \, dt =
\left( \operatorname{ sign } \lambda_j(A) \right) \infty,
\ \ \inf\limits_{ t \in I } \left| \hat p_{ jj }(t) \right| > 0
\ \ \left( j = \overline{ 1, n } \right).
\]
В силу формул \eqref{NPI_q(t)}, \eqref{NPI_P(t)}, \eqref{IE:NPI_Lip} функции
\[
\hat q_j(t) = O \left( 1 \right),
\ \ \hat p_{ jl }(t) = O(1),
\ \ \hat x(t) = O \left( \left( \varepsilon_1(t)
\right)^{ 1 + \frac{ \rho_{ k+1 } }{ \varrho_1 } } \right),
\]
поэтому для достаточно малого числа $ d_0 $ отношения
\[
\dfrac{ \hat q_j( t ) }{ \hat p_{ jj }(t) } = O \left( 1 \right),
\ \ \dfrac{ \hat p_{ jl }(t) }{ \hat p_{ jj }(t) } = O(1),
\ \ \dfrac{ \hat x( t ) }{ \hat p_{ jj }(t) } = o(1)
\ \ \left( t \to + \infty; \ j, l = \overline{ 1, n }; \ j \ne l \right).
\]
Значит, для выражений $ \hat q_j^*(t) $, $ \hat p_{ jl }^*(t) $,
$ \hat x_j^*(t) $ $ \left( j, l = \overline{ 1, n }; \ j \ne l \right) $
выполнены условия лемм 2 и 2$'$ из \S \,~2 главы II кандидатской диссертации
\cite{Bib:Kostin_PHDT} (страницы 45 и 47 соответственно).
Следовательно, для достаточно малого числа $ d_0 $ выражения
$ \hat q_j^*(t) = O(1), $ $ \hat p_{ jl }^*(t) = O(1), $
$ \hat x_j^*(t) = o(1) $
$ \left( t \to + \infty; \ j, l = \overline{ 1, n }; \ j \ne l \right). $
Так как подынтегральные функции во внешних интегралах в этих выражениях
неотрицательны, то за счёт выбора числа $ t_0 $ достаточно большим, а $ d_0 $
достаточно малым можем сделать величины
\begin{gather*}
\hat q_0 := \max\limits_{ j = \overline{ 1, n } }
\left( \sup\limits_{ t \in I } \hat q_j^*(t) \right) \ge 0,
\ \ \hat p_1 := \max\limits_{ \substack{ j,l = \overline{ 1, n } \\ j > l } }
\left( \sup\limits_{ t \in I } \hat p_{ jl }^*(t) \right) \ge 0,
                                                                        \\
\hat p_2 := \max\limits_{ \substack{ j,l = \overline{ 1, n } \\ j < l } }
\left( \sup\limits_{ t \in I } \hat p_{ jl }^*(t) \right) \ge 0,
\ \ \hat x_0 := \max\limits_{ j = \overline{ 1, n } }
\left( \sup\limits_{ t \in I } \hat x_j^*(t) \right) \ge 0
\end{gather*}
сколь угодно малыми.
Очевидно, что пределы $ \lim\limits_{ t \to + \infty } \hat z_j^*(t) = 0 $
$ \left( j = \overline{ 1, n } \right), $ поэтому для любых начальных
значений $ z_j(t_0) \in \mathbb R $ (индексы $ j $ такие, что справедливы
требования \eqref{C:NPI_No_of_parameters_O}) найдется такое достаточно
большое число $ t_1 \in \mathbb R $ ($ t_1 \ge t_0 $), что величина
\[
\hat z_0 := \max\limits_{ j = \overline{ 1, n } }
\left( \sup\limits_{ t \ge t_1 } \hat z_j^*(t) \right) \ge 0
\]
тоже будет сколь угодно малой.

Оценим при $ t \ge t_1 $ функции $ \hat \xi_j(t, t_0, \hat \epsilon_0) $
$ \left( j = \overline{ 1, n } \right). $
Нетрудно видеть, что при $ t \ge t_1 $ функции
$ \hat \xi_j(t, t_0, \hat \epsilon_0) $
$ \left( j = \overline{ 1, n } \right) $ мажорируются константами
$ \hat \xi_j $ $ \left( j = \overline{ 1, n } \right) $, которые определяются
таким образом:
\begin{align*}
& \hat \xi_j = \hat q_0 + ( n-1 ) \hat p_1 \hat \epsilon_0 + \hat z_0 +
\hat x_0 \hat \epsilon_0 + \hat p_2 \sum_{ l = j+1 }^n \hat \xi_l
\ \ \left( j = \overline{ 1, n-1 } \right),
                                                                        \\
& \hat \xi_n = \hat q_0 + ( n-1 ) \hat p_1 \hat \epsilon_0 + \hat z_0 +
\hat x_0 \hat \epsilon_0.
\end{align*}
Решая эти уравнения, получим
\[
\hat \xi_j = ( 1 + \hat p_2 )^{ n-j } \hat \xi_n
\ \ \left( j = \overline{ 1, n } \right).
\]
Так как константы $ \hat \xi_n $ и $ \hat p_2 $ здесь неотрицательны, то
ясно, что среди $ \hat \xi_j $ $ \left( j = \overline{ 1, n } \right) $
наибольшую величину будет иметь $ \hat \xi_1. $
Поэтому условие ограниченности при $ t \ge t_1 $ последовательных приближений
одной и той же константой $ \hat \epsilon_0 $ будет иметь следующий вид:
\[
\hat \xi_1 = ( 1 + \hat p_2 )^{ n-1 } \left( \hat q_0 +
( n-1 ) \hat p_1 \hat \epsilon_0 + \hat z_0 + \hat x_0 \hat \epsilon_0
\right) \le \hat \epsilon_0
\]
или
\[
\hat z_0^* := \dfrac{ \hat q_0 + \hat z_0 }{ \dfrac{ 1 }{ ( 1 + \hat p_2
)^{ n-1} } - ( n-1 ) \hat p_1 - \hat x_0 } \le \hat \epsilon_0.
\]
Учитывая, что для достаточно малой постоянной $ d_0 $ и достаточно больших
чисел $ t_0, $ $ t_1 $ величины $ \hat q_0, $ $\hat p_1, $ $ \hat p_2, $
$ \hat x_0, $ $ \hat z_0 $ сколь угодно малы, можем сделать число
$ \hat z_0^* $ не превосходящим константы $ \hat a $ (где $ \hat a $
определяет область $ \hat D_3 $).
Поэтому в качестве $ \hat \epsilon_0 $ можно взять число $ \hat z_0^* $.
Так как для достаточно малой постоянной $ d_0 $ и достаточно большого числа
$ t_0 $ величины $ \hat p_1, \ \hat x_0 $ сколь угодно малы, то выполнимость
условия
\[
\max\limits_{ j = \overline{ 1, n } } \left( \sup\limits_{ t \in I } \left(
\sum_{ l=1 }^{ j-1 } \hat p_{ jl }^*(t) + \hat x_j^*(t) \right) \right) < 1
\]
в данном случае очевидна.
Таким образом, для системы \eqref{NPI_xi_j(t t_0 epsilon_0)} выполнены все
условия теоремы 1.1 из \S \ 1 главы IV кандидатской диссертации
\cite{Bib:Kostin_PHDT} (страница 67).
Поэтому, для достаточно большого числа $ t_0 $ система дифференциальных
уравнений \eqref{SDE:NPI_solution->0} будет заведомо иметь хотя бы одно
вещественное ограниченное на промежутке $ I $ частное решение $ z_j(t) $
$ \left( j = \overline{ 1, n } \right) $ с условием
$ | z_j(t) | \le \hat z_0^* $
$ \left( t \in I; \ j = \overline{ 1, n } \right) $.
Более того, на промежутке $ [ t_1, + \infty ) $ решение $ z_j(t) $
$ \left( j = \overline{ 1, n } \right) $ зависит от стольких произвольных
постоянных, сколько имеется индексов $ j $ таких, что выполнены условия
\eqref{C:NPI_No_of_parameters_O} ($ t_1 \in \mathbb R $ --- достаточно
большое число, определяется этими постоянными, $ t_1 \ge t_0 $).
Возвращаясь обратно к вектор--функции $ r $, получим требуемое.

Как известно, общее решение системы обыкновенных дифференциальных уравнений
$ n $--го порядка зависит от $ n $ произвольных параметров, следовательно в
семействе решений вида \eqref{CV:r_j} могут фигурировать не более чем $ n $
независимых скалярных постоянных.

\end{proof}

\subsection{Случай простых собственных значений матрицы A}

В этой части статьи исследован асимптотический характер формального частного
решения (семейства решений) \eqref{FSolution} основного вещественного
векторного квазилинейного обыкновенного дифференциального уравнения
\eqref{DEV:main} в том случае, когда матрица $ A $ не имеет кратных
собственных значений (теоремы \ref{T:AsymptChart_O} -- \ref{T:AsymptChart}).

\subsubsection{Вспомогательная лемма о выделении суммируемых слагаемых}

Ниже используем известные обозначения.
Пусть число $ p \in \mathbb R^+ $ и
\begin{gather*}
\mathcal R(I) := \left\{ g(t) : I \to \mathbb C \ \left|
\ \exists \int\limits_I g(t) \, d t \in \mathbb C \right. \right\},
                                                                        \\
\mathcal R_p(I) := \left\{ g(t) : I \to \mathbb C \ \Big|
\ \left| g(t) \right|^p \in \mathcal R(I) \ \right\}.
\end{gather*}

В следующей лемме выделены заведомо суммируемые на промежутке $ I $ слагаемые
в линейной однородной части векторного дифференциального уравнения
\eqref{DEV:separatedFT}.

\begin{lemma}
                                                    \label{L:ready>L-Diag}
Пусть для векторного дифференциального уравнения \eqref{DEV:main} выполнены
все условия леммы \ref{L:separatedFT}, а также справедливы следующие
предположения:
\begin{enumerate}
\item
существует число $ \varpi \in \mathbb R^+ $ такое, что функция
$ \varepsilon_1(t) \in \mathcal R_\varpi(I) $;
\item
производные $ \dfrac{ d^r \varepsilon_l(t) }{ dt^r } \in \mathcal R_1(I) $
($ r = \overline{ 1, \left[ \rho_k - \varrho_l \right] } $,
$ \rho_k - \varrho_l \ge 1 $);
\item
вектор--функции
$ { f }_l( t, { y } ) \in C_{ t, \ { y } }^{ 0, \ s_l + 2 } ( D ) $,
числа $ s_l := \left[ \varpi - \frac { \varrho_l } { \varrho_1 } \right] $
($ \varrho_l \le \varpi \varrho_1 $);
\item
супремумы
$ \sup\limits_{ ( t,  y ) \in D } \left\|
\partial_{ { y } }^\wp { f }_l( t, { y } ) \right\| < + \infty $
$ \left( \wp \in \mathbb N_0^n, \ | \wp | = \overline{ 3, s_l + 2 },
\ \varrho_l \le \varpi \varrho_1 \right) $;
\item
вектор--функции
$ \partial_{ { y } }^\wp { f }_l( t, \varphi_0( t, c_0 ) ) \in
\mathcal K_1^{ n \times 1 }(A) $
$ \left( \wp \in \mathbb N_0^n, \ | \wp | = \overline{ 1, s_l + 1 },
\ \varrho_l \le \varpi \varrho_1 \right) $.
\end{enumerate}
Тогда для достаточно большого числа $ t_0 $ с помощью замены неизвестной
вектор--функции вида \eqref{CV:r_j}, вводя новую неизвестную вектор--функцию
$ r : I \to \mathbb R^{ n \times 1 } $, векторное дифференциальное уравнение
\eqref{DEV:main} можно привести к эквивалентному дифференциальному уравнению
вида
\begin{gather}
\dfrac{ dr }{ dt } = \left( A + A( t ) + \mathfrak S(t) \right) r +
 g(t) + \sum_{ l=1 }^m \varepsilon_l(t) { \psi }_l( t, { r } )
\ \ \left( \left( t, { r } \right) \in D_1 \right),
                                                    \label{DEV:preparedLHP}
                                                                        \\
A( t ) = \sum_{ \substack{ 1 \le | \wp | < \varpi \\
\wp \in \mathbb N_0^m \\ ( \mu_\wp(t) \notin \mathcal R_1(I) ) } }
\mu_\wp(t) A_\wp( t ),
\ \ \mu_\wp(t) = \prod_{ l=1 }^m
\left( \varepsilon_l(t) \right)^{ ( \wp )_l },
                                                    \label{Delta(t)}
\end{gather}
где матрицы $ A_\wp( t ) \in \mathcal K_1^{ n \times n }(A) $
($ 1 \le | \wp | < \varpi, $ $ \wp \in \mathbb N_0^m $,
$ \mu_\wp(t) \notin \mathcal R_1(I) $),
$ \mathfrak S(t) \in \mathcal R_1^{ n \times n }(I) \cap C(I). $

\end{lemma}

\begin{proof}

В силу леммы \ref{L:separatedFT} для достаточно большого числа $ t_0 $ с
помощью замены неизвестной вектор--функции вида \eqref{CV:r_j}, вводя новую
неизвестную вектор--функцию $ { r } : I \to \mathbb R^{ n \times 1 } $,
векторное дифференциальное уравнение \eqref{DEV:main} можно привести к
эквивалентному дифференциальному уравнению вида \eqref{DEV:separatedFT}.
Рассмотрим линейную однородную часть векторного дифференциального уравнения
\eqref{DEV:separatedFT}.
Выделим из суммы
\begin{equation}
\sum_{ l=1 }^m \varepsilon_l(t)
\frac { \partial { f }_l ( t, s( t ) ) } { \partial { y } }
                                                        \label{E:AltTerm}
\end{equation}
заведомо суммируемые на промежутке $ I $ слагаемые.
Для этого нам потребуются нижеследующие свойства.

Заметим, что функции
\begin{equation}
\nu_{ sp }(t) \in \mathcal R_1(I) \ \ (\rho_s \ge \varpi \varrho_1,
\ p = \overline{ 1, \varkappa_s }).
                                                        \label{P:nu_in_R1}
\end{equation}
Действительно, в силу оценок \eqref{estimate_nu} имеем свойства
$ \nu_{ sp }(t) =
O \left( \left( \varepsilon_1(t) \right)^{ \varpi } \right) $
($ t \to + \infty $, $ \rho_s \ge \varpi \varrho_1 $,
$ p = \overline{ 1, \varkappa_s } $).
Следовательно, функции $ | \nu_{ sp }(t) | $ ($ \rho_s \ge \varpi \varrho_1 $,
$ p = \overline{ 1, \varkappa_s } $) для достаточно больших $ t \in I $
ограничены функцией из класса $ \mathcal R_1(I) $.

В силу свойств \eqref{P:nu_in_R1} слагаемые из суммы \eqref{E:AltTerm},
содержащие функции $ \varepsilon_l(t) $, ранга
$ \varrho_l \ge \varpi \varrho_1 $, будут суммируемыми на промежутке $ I $,
так как
\begin{equation*}
C(I) \ni \frac { \partial { f }_l ( t, s( t ) ) }
{ \partial { y } } = O \left( 1 \right)
\ \ ( t \to + \infty, \ l = \overline{ 1, m } ).
\end{equation*}

Рассмотрим случай, когда $ \varrho_l < \varpi \varrho_1 $.
Учитывая оценки \eqref{estimate_delta_K}, с помощью формулы Тейлора получим
разложения матриц Якоби
\begin{equation*}
\frac { \partial { f }_l ( t, s( t ) ) } { \partial { y } } =
\sum_{ \substack{ 0 \le | \wp | \le s_l \\ ( \wp \in \mathbb N_0^n ) } }
\frac 1 { \wp ! } \frac { \partial } { \partial { y } }
\partial_{ { y } }^\wp { f }_l ( t, \varphi_0( t, c_0 ) ) \delta^\wp( t ) +
\left( \varepsilon_1(t) \right)^{ s_l + 1 } O \left( 1 \right)
\ \ ( t \to + \infty, \ \varrho_l < \varpi \varrho_1 ).
\end{equation*}
Подставляя эти разложения в сумму \eqref{E:AltTerm}, получим матричное
равенство
\begin{multline}
\sum_{ \varrho_l < \varpi \varrho_1 } \varepsilon_l(t)
\frac { \partial { f }_l ( t, s( t ) ) } { \partial { y } } =
\sum_{ \varrho_l < \varpi \varrho_1 } \varepsilon_l(t)
\sum_{ \substack{ 0 \le | \wp | \le s_l \\ ( \wp \in \mathbb N_0^n ) } }
\frac 1 { \wp ! } \frac { \partial } { \partial { y } }
\partial_{ { y } }^\wp { f }_l ( t, \varphi_0( t, c_0 ) ) \delta^\wp( t ) +
                                                                        \\
+ \sum_{ \varrho_l < \varpi \varrho_1 } \varepsilon_l(t)
\left( \varepsilon_1(t) \right)^{ s_l + 1 } O \left( 1 \right)
\ \ ( t \to + \infty ).
                                                    \label{AltTermPart2}
\end{multline}
Оценим ранг
\begin{multline*}
\operatorname{ Rank } \left( \varepsilon_l(t)
\left( \varepsilon_1(t) \right)^{ s_l + 1 } \right) =
\varrho_l + \varrho_1 \left( s_l + 1 \right) =
\varrho_l + \varrho_1 \left( \left[ \varpi -
\frac { \varrho_l } { \varrho_1 } \right] + 1 \right) =
                                                                        \\
= \varrho_l + \varrho_1 \left( \varpi -
\frac { \varrho_l } { \varrho_1 } - \left\{ \varpi -
\frac { \varrho_l } { \varrho_1 } \right\} + 1 \right) =
\varpi \varrho_1 + \varrho_1 \left( 1 - \left\{ \varpi -
\frac { \varrho_l } { \varrho_1 } \right\} \right) >
\varpi \varrho_1
\ \ ( \varrho_l < \varpi \varrho_1 ).
\end{multline*}
Поэтому, принимая во внимание свойства \eqref{P:nu_in_R1}, заметим что сумма,
содержащая символ $ O(1) $ в правой части матричного равенства
\eqref{AltTermPart2} будет из класса
$ \mathcal R_1^{ n \times n }(I) \cap C(I) $.

Рассмотрим двойную сумму в правой части матричного равенства
\eqref{AltTermPart2}.
Учитывая аксиомы класса $ \mathcal K_1(A) $ можем заключить, что произведение
двух любых компонент вектор--функции $ \delta(t) $ можно записать в виде
суммы, такого же вида, как и у этих множителей.
Поэтому в силу условий леммы столбцы этой двойной суммы обладают формой такой
же, как и у вектор--функции $ \delta(t) $.
Таким образом, эту двойную сумму можно переписать в следующем виде
\begin{equation}
\sum_{ s=1 }^{ \varsigma } \sum_{ p=1 }^{ \varkappa_s }
\nu_{ sp }(t) A_{ sp }( t ),
                                                    \label{AltTerm_gamma}
\end{equation}
где $ \varsigma \in \mathbb N $ --- некоторое известное число; в
представлениях \eqref{E:nu_sp(t)}, функций $ \nu_{ sp }(t) $
($ s = \overline{ 1, \varsigma } $, $ p = \overline{ 1, \varkappa_s } $),
могут участвовать производные функций $ \varepsilon_l(t) $
($ \rho_k - \varrho_l \ge 1 $) порядка не выше чем
$ \left[ \rho_k - \varrho_l \right] $;
$ A_{ sp }( t ) \in \mathcal K_1^{ n \times n }(A) $
($ s = \overline{ 1, \varsigma } $, $ p = \overline{ 1, \varkappa_s } $) ---
известные матрицы.

Разобьем сумму \eqref{AltTerm_gamma} на несуммируемые и суммируемые на
промежутке $ I $ составляющие.
Во первых, отнесем ко второй составляющей слагаемые из выражения
\eqref{AltTerm_gamma}, содержащие в качестве коэффициентов функции
$ \nu_{ sp }(t) $, ранга $ \rho_s \ge \varpi \varrho_1 $
($p = \overline{ 1, \varkappa_s }$), так как они будут суммируемыми на
промежутке $ I $ в силу свойств \eqref{P:nu_in_R1} и условий леммы.
Во вторых, отнесем туда же все члены из суммы \eqref{AltTerm_gamma}
содержащие в представлениях \eqref{E:nu_sp(t)} функций $ \nu_{ sp }(t) $
($ s = \overline{ 1, \varsigma } $, $ p = \overline{ 1, \varkappa_s } $)
производные функций $ \varepsilon_l(t) $ ($ \rho_k - \varrho_l \ge 1 $), так
как по условию леммы они будут суммируемыми на промежутке $ I $.
В третьих, в выражении \eqref{AltTerm_gamma} могут быть суммируемые на
промежутке $ I $ функции $ \nu_{ sp }(t) $, ранга
$ \rho_s < \varpi \varrho_1 $ ($ p \in \overline{ 1, \varkappa_s } $) и
несодержащие производных функций $ \varepsilon_l(t) $
($ l = \overline{ 1, m } $).

Оставшиеся после этого слагаемые в сумме \eqref{AltTerm_gamma} будут
содержать функции $ \nu_{ sp }(t) \notin \mathcal R_1(I) $ (ранга
$ \rho_s < \varpi \varrho_1 $), состоящие только из произведений некоторых
степеней (обозначим их как $ \beta_{ spl } \in \mathbb N_0 $) функций
$ \varepsilon_l(t) \notin \mathcal R_1(I) $.
Рассмотрим ранг функций $ \nu_{ sp }(t) $ в этих слагаемых:
$$
\sum \varrho_l \beta_{ spl } < \varpi \varrho_1
\ \left( \nu_{ sp }(t) \notin \mathcal R_1(I) \right).
$$
Разделим обе части этого неравенства на число $ \varrho_1 > 0 $ и оценим
снизу его левую часть, учитывая, что
$ 1 \le \frac{ \varrho_l }{ \varrho_1 } $ ($ l = \overline{ 1, m } $),
получим
$$
\sum \beta_{ spl } < \varpi
\ \left( \nu_{ sp }(t) \notin \mathcal R_1(I) \right).
$$

Перенумеруем матрицы $ A_{ sp }( t ) $ в несуммируемых на промежутке
$ I $ членах из выражения \eqref{AltTerm_gamma} с помощью мультииндексов из
множества $ \mathbb N_0^m $ в соответствии со степенями функций
$ \varepsilon_l(t) $ ($ l = \overline{ 1, m } $) стоящими перед этими
матрицами.
В силу последних неравенств модули этих мультииндексов будут меньше числа
$ \varpi $.

Таким образом доказали, что сумму \eqref{AltTerm_gamma} можно представить в
виде $ A( t ) + \mathfrak S_1 ( t ), $ где $ A( t ) $ определяется по формуле
\eqref{Delta(t)},
$ \mathfrak S_1 ( t ) \in \mathcal R_1^{ n \times n }(I) \cap C(I) $
--- известные матрицы.
Итак, векторное дифференциальное уравнение \eqref{DEV:separatedFT} имеем
право переписать в виде \eqref{DEV:preparedLHP}.

\end{proof}

\begin{remark}

Сопоставляя двойную сумму из правой части матричного равенства
\eqref{AltTermPart2} с выражением \eqref{AltTerm_gamma} можно легко найти
формулы для матриц $ A_\wp( t ) $ в первых слагаемых суммы $ A( t ) $:
\begin{equation}
A_\wp( t ) = \frac { \partial } { \partial { y } }
{ f }_l ( t, \varphi_0( t, c_0 ) ),
\ \wp = ( \underbrace{ 0, \dots, 0 }_{ l-1 }, 1, 0, \dots, 0 )
\ \left( \wp \in \mathbb N_0^m, \ \varepsilon_l(t) \notin \mathcal R_1(I)
\right).
                                                    \label{A_p(t), |p|=1}
\end{equation}

\end{remark}

\subsubsection{\texorpdfstring{Абстрактные классы $ \mathcal K_2(A) $ и
$ \mathcal K_3(A) $ колеблющихся функций}{Второй и третий абстрактные классы
колеблющихся функций}}

Важную роль играют матрица и преобразование Ляпунова (смотрите
\cite{Bib:Demidovich}).

\begin{definition}

Матрица $ L(t) : I \to \mathbb C^{ n \times n } $, $ L(t) \in C^1(I) $
называется матрицей Ляпунова, если выполнены следующие условия:
$$
\sup\limits_{ t \in I } || L(t) || < + \infty,
\ \sup\limits_{ t \in I } \left\| \dfrac{ d L(t) }{ dt } \right\| < + \infty,
\ | \det L(t) | \ge const > 0 \ ( \forall t \in I ).
$$

\end{definition}

\begin{definition}

Линейное преобразование $  y = L(t)  x $ с матрицей Ляпунова
$ L(t) : I \to \mathbb C^{ n \times n } $, где вектор--функции
$  x,  y : I \to \mathbb C^{ n \times 1 } $,
называется преобразованием Ляпунова.

\end{definition}

Аксиоматически введем в рассмотрение абстрактный класс $ \mathcal K_2(A) $
колеблющихся функций, который необходим, для того чтобы асимптотически
привести линейную однородную часть векторного дифференциального уравнения
типа \eqref{DEV:preparedLHP} (с условием $ 0 \notin \Delta( A ) $) с
помощью вещественного преобразования Ляпунова $ Q(t) $ (более того
$ \exists \lim\limits_{ t \to + \infty } \det Q(t) \ne 0 $) к диагональному
виду такому же, как и у соответствующей части уравнения
\eqref{DEV:DiagLH_Part} (смотрите \cite{Bib:Odessa01}).

\begin{definition}

Будем обозначать через $ \mathcal K_2(A) $ некоторый непустой класс
функций $ \left\{ f(t) : I \to \mathbb R \right\} $ таких, что
\begin{enumerate}
\item
выполнены свойства \ref{C:1_aca}. и \ref{C:2_aca}. из
определения \ref{D:aca} с заменой класса $ \mathcal K_1(A) $ на
$ \mathcal K_2(A) $;
\item
если функция $ f(t) \in \mathcal K_2(A) $, то существует её среднее значение
$ \operatorname{ M } \left( f(t) \right) $ и
$$
\sup\limits_{ t \in I } \left| \int\limits_{ t_0 }^t
\left( f( \tau ) - \operatorname{ M } \left( f( \tau ) \right) \right)
\, d \tau \right| < + \infty;
$$
\item
выполнено свойство \ref{C:3_aca}. из определения \ref{D:aca} с заменой класса
$ \mathcal K_1(A) $ на $ \mathcal K_2(A) $ и множества $ \Lambda( A ) $ на
 $ \Delta( A ) $.
\end{enumerate}

\end{definition}

Приведем три важных примера классов $ \mathcal K_2(A) $ аналогичных ранее
указанным примерам классов $ \mathcal K_1(A) $.

\begin{example}
                                                    \label{Ex:PeriodicK2}
Заменяя всюду в примере \ref{Ex:PeriodicK1} класс $ \mathcal K_1(A) $ на
$ \mathcal K_2(A) $ и множество $ \Lambda( A ) $ на $ \Delta( A ) $,
получим пример периодического класса $ \mathcal K_2(A) $.

\end{example}

\begin{example}
                                                \label{Ex:AlmostPeriodicK2}
Проделав тоже самое (смотрите предыдущий пример), но в примере
\ref{Ex:AlmostPeriodicK1}, а также добавив ограничение на множество
$ \Gamma_2 $
\begin{equation*}
\inf_{ \gamma \in \Gamma_2 \setminus \{ 0 \} }
\left| \gamma \right| > 0,
\end{equation*}
получим пример РПП класса $ \mathcal K_2(A) $.

\end{example}

\begin{example}
                                            \label{Ex:Decreasing_Exp_LC_K2}
Проделав тоже самое с заменой множества $ \Gamma_2 $ на $ \Gamma_3 $ (смотрите
предыдущий пример), но в примере \ref{Ex:Decreasing_Exp_LC_K1}, получим
пример класса $ \mathcal K_2(A) $, состоящего из линейных комбинаций
экспоненциальных функций.

\end{example}

Аксиоматически введем в рассмотрение абстрактный класс $ \mathcal K_3(A) $
колеблющихся функций, который необходим, для того чтобы построить, по крайней
мере, одно формальное частное решение вида \eqref{FSolution} основного
векторного дифференциального уравнения \eqref{DEV:main} и асимптотически
привести к вещественному диагональному виду линейную однородную часть
векторного дифференциального уравнения типа \eqref{DEV:preparedLHP} (с
условием $ 0 \notin \Delta( A ) $).

\begin{definition}

Будем обозначать через $ \mathcal K_3(A) $ некоторый непустой класс
функций, удовлетворяющий аксиоматике классов $ \mathcal K_1(A) $ и
$ \mathcal K_2(A) $ одновременно (причём во всех аксиомах этих классов
классы $ \mathcal K_1(A) $ и $ \mathcal K_2(A) $ заменяются на
$ \mathcal K_3(A) $).

\end{definition}

Как и прежде приведем три важных примера классов $ \mathcal K_3(A) $
аналогичных ранее указанным примерам классов $ \mathcal K_1(A) $.

\begin{example}
                                                    \label{Ex:PeriodicK3}
Заменяя всюду в примере \ref{Ex:PeriodicK1} класс $ \mathcal K_1(A) $ на
$ \mathcal K_3(A) $ и множество $ \Lambda( A ) $ на
$ ( \Lambda( A ) \cup \Delta( A ) ) $, получим пример периодического
класса $ \mathcal K_3(A) $.

\end{example}

\begin{example}
                                                \label{Ex:AlmostPeriodicK3}
Проделав тоже самое (смотрите предыдущий пример), но в примере
\ref{Ex:AlmostPeriodicK1}, а также добавив ограничение на множество
$ \Gamma_2 $
\begin{equation*}
\inf_{ \gamma \in \Gamma_2 \setminus \{ 0 \} }
\left| \gamma \right| > 0,
\end{equation*}
получим пример РПП класса $ \mathcal K_3(A) $.

\end{example}

\begin{example}
                                            \label{Ex:Decreasing_Exp_LC_K3}
Проделав тоже самое с заменой множества $ \Gamma_2 $ на $ \Gamma_3 $ (смотрите
предыдущий пример), но в примере \ref{Ex:Decreasing_Exp_LC_K1}, получим
пример класса $ \mathcal K_3(A) $, состоящего из линейных комбинаций
экспоненциальных функций.

\end{example}

\begin{remark}

Очевидно, что указанными примерами классов $ \mathcal K_1(A) $,
$ \mathcal K_2(A) $ и $ \mathcal K_3(A) $ не исчерпываются все возможные их
варианты.
Действительно, при некоторых условиях можно привести примеры этих классов
состоящих из конечных линейных комбинаций функций вида
$ e^{ \gamma_1 t } \sin \alpha_1 t $, $ e^{ \gamma_2 t } \cos \alpha_2 t $
(числа $ \gamma_s \in \mathbb R_0^- ( \mathbb R^- ) $,
$ \alpha_s \in \mathbb R $, $ s = \overline{ 1, 2 } $) с вещественными
коэффициентами.
Тем самым применённый в работе аксиоматический подход позволяет охватить
достаточно широкий набор функций.

\end{remark}

\subsubsection{Теоремы об асимптотическом характере}

В следующей теореме исследован асимптотический характер формального частного
решения (семейства решений) \eqref{FSolution} основного вещественного
векторного квазилинейного обыкновенного дифференциального уравнения
\eqref{DEV:main} в том случае, когда матрица $ A $ не имеет кратных
собственных значений, а оценка погрешности $ r $ получена в виде $ O $.

\begin{theorem}

Пусть для векторного дифференциального уравнения \eqref{DEV:main} выполнены
все условия леммы \ref{L:ready>L-Diag} за исключением того, что в них всюду
класс $ \mathcal K_1(A) $ заменяется на $ \mathcal K_3(A) $, а также
справедливы следующие предположения:
                                                    \label{T:AsymptChart_O}
\begin{enumerate}
\item
число $ 0 \notin \Delta( A ) $;
\item
функции $ \varepsilon_l(t) \downarrow \uparrow 0 $,
$ \dfrac{ d \varepsilon_l(t) }{ dt } \in C(I) \cap \mathcal R_1(I) $,
$ \varepsilon_1(t) > 0 $ ($ \varepsilon_l(t) \notin \mathcal R_1(I) $,
$ t \in I $, $ t \to + \infty $);
\item
любая функция вида
\[
\gamma_1 \left( \varepsilon_1(t) \right)^{ -1 }
\dfrac{ d \varepsilon_1(t) }{ dt } +
\sum_{ \substack{ 1 \le | \wp | < \varpi \\ \wp \in \mathbb N_0^m
\\ ( \mu_\wp(t) \notin \mathcal R_1(I) ) } }
\gamma_\wp \mu_\wp( t ),
\]
где числа $ \gamma_1 $, $ \gamma_\wp \in \mathbb R $
$ ( 1 \le | \wp | < \varpi $, $ \wp \in \mathbb N_0^m $,
$ \mu_\wp(t) \notin \mathcal R_1(I) ) $,
сохраняет знак, по крайней мере, в нестрогом смысле, для достаточно больших
значений аргумента $ t \in I $;
\item
число $ k $ такое, что выполнены неравенства
$ \rho_{ k+1 } > \varrho_1 \varpi $,
$ \rho_{ k+1 } \ge \varrho_1 ( 2 \varpi - 1 ) $.
\end{enumerate}
Тогда для достаточно большого числа $ t_0 $ у векторного дифференциального
уравнения \eqref{DEV:main} на промежутке $ I $ существует хотя бы одно
частное решение вида \eqref{CV:r_j}, где вектор--функции
$ \varphi_0( t, c_0 ), \varphi_{ sp }( t, c_{ sp } ) \in
\mathcal K_3^{ n \times 1 }(A) \cap C^1( I ) $
($ s = \overline{ 1, k } $, $ p = \overline{ 1, \varkappa_s } $), причём
погрешность
\begin{equation}
r = O \left( \left( \varepsilon_1(t) \right)^{
\frac { \rho_{ k+1 } } { \varrho_1 } - \varpi } \right)
\ \ ( t \to + \infty ).
                                                \label{r_j = O}
\end{equation}
Более того, погрешность $ r $ зависит от стольких достаточно малых по
абсолютной величине произвольных скалярных постоянных, сколько имеется
индексов $ j $ таких, что справедливы условия
\begin{equation}
j \in \left\{ 1, 2, 3, \dots, n \right\},
\ \ \int\limits_I p_{ jj }(t) \, dt \ne + \infty,
                                                \label{No_of_parameters_O}
\end{equation}
где функции $ p_{ jj }(t) $ ($ j = \overline{ 1, n } $) определяются по
формулам \eqref{p_jj(t)}.
(В семействе решений вида \eqref{CV:r_j} могут фигурировать не более чем
$ n $ независимых скалярных параметров.)

\end{theorem}

\begin{proof}

Заметим, что леммы \ref{L:truncated_main}, \ref{L:separatedFT} и
\ref{L:ready>L-Diag}, теорема \ref{T:FSolution} остаются в силе, если в них
всюду класс $ \mathcal K_1(A) $ заменить на $ \mathcal K_3(A) $.
Переформулируем их для класса $ \mathcal K_3(A) $ вместо $ \mathcal K_1(A) $.
В силу леммы \ref{L:ready>L-Diag} для класса $ \mathcal K_3(A) $, для
достаточно большого числа $ t_0 $ с помощью замены неизвестной
вектор--функции вида \eqref{CV:r_j} (причём в сумме \eqref{s_j}
вектор--функции $ \varphi_0( t, c_0 ) $, $ \varphi_{ sp }( t, c_{ sp } ) $
$ \left( s = \overline{ 1, k }, \ p = \overline{ 1, \varkappa_s } \right) $
нужно взять из класса $ \mathcal K_3^{ n \times 1 }(A) $), вводя новую
неизвестную вектор--функцию $ r : I \to \mathbb R^{ n \times 1 } $, векторное
дифференциальное уравнение \eqref{DEV:main} можно привести к эквивалентному
векторному дифференциальному уравнению вида \eqref{DEV:preparedLHP}.
Асимптотически приведем линейную однородную часть этого векторного
дифференциального уравнения к диагональному виду с помощью теоремы 2 из
работы авторов \cite{Bib:Odessa01} (смотрите также \cite{Bib:RNASU04}).
Убедимся, что выполнены все её условия.

Докажем, что для любого мультииндекса $ \wp \in \mathbb N_0^m $,
$ | \wp | \ge \varpi $ любая функция $ \mu_\wp (t) $ будет суммируемой на
промежутке $ I $.
Для этого оценим снизу ранг функций $ \mu_\wp (t) $
($ \wp \in \mathbb N_0^m $, $ | \wp | \ge \varpi $), учитывая, что
$ \frac{ \varrho_l }{ \varrho_1 } \ge 1 $ ($ l = \overline{ 1, m } $),
получим
$$
\operatorname{ Rank } ( \mu_\wp (t) ) = \sum_{ l=1 }^m \varrho_l ( \wp )_l =
\varrho_1 \sum_{ l=1 }^m \frac{ \varrho_l }{ \varrho_1 } ( \wp )_l \ge
\varrho_1 \sum_{ l=1 }^m ( \wp )_l =
\varrho_1 | \wp | \ge \varrho_1 \varpi.
$$
Таким образом, в силу свойств \eqref{P:nu_in_R1} получили требуемое.

И так, выполнены все требования упомянутой теоремы.
Сделаем в векторном дифференциальном уравнении \eqref{DEV:preparedLHP}
замену неизвестной вектор--функции вида (21) (с матрицей перехода Ляпунова
$ Q(t) : I \to \mathbb R^{ n \times n } $,
$ \exists \lim\limits_{ t \to + \infty } \det Q(t) \ne 0 $) из теоремы 2
вышеуказанной публикации авторов.
В результате получим векторное дифференциальное уравнение вида
\begin{equation}
\dfrac{ dx }{ dt } = \Re \, D(t) x + Q^{ -1 }(t) g(t) +
Q^{ -1 }(t) \sum_{ l=1 }^m \varepsilon_l(t) \psi_l( t, Q(t) x )
\ \ \left( \left( t, x \right) \in D_2 \right),
                                                \label{DEV:DiagLH_Part}
\end{equation}
где матрица $ D(t) $ определяется по формуле (6) из статьи авторов
\cite{Bib:Odessa01}:
\begin{equation*}
D(t) := \operatorname{ diag } ( d_1(t), \dots, d_n(t) ),
\ \ d_j(t) := \lambda_j(A) +
\sum_{ \substack{ 1 \le | \wp | < \varpi \\
\wp \in \mathbb N_0^m \\ ( \mu_\wp(t) \notin \mathcal R_1(I) ) } }
\gamma_{ j \wp } \mu_\wp( t ) + \varsigma_j(t),
\end{equation*}
где числа $ \gamma_{ j \wp } \in \mathbb C $, функции
$ \varsigma_j(t) : I \to \mathbb C $,
$ \varsigma_j(t) \in \mathcal R(I) \cap C(I) $ ($ j, \wp $ любые допустимые),
область
$ D_2 := \left\{ ( t, x ) \ \left| \ t \in I,
\ x \in \mathbb R^{ n \times 1 }, \ || x || \le \hat a_0 \left( n
\sup\limits_{ t \in I } || Q(t) || \right)^{ -1 } \right. \right\}. $

Добьемся суммируемости на промежутке $ I $ свободного члена этого векторного
дифференциального уравнения.
Для этого сделаем в векторном дифференциальном уравнении
\eqref{DEV:DiagLH_Part} замену неизвестной вектор--функции вида
\[
x = \left( \varepsilon_1(t) \right)^{
\frac { \rho_{ k+1 } } { \varrho_1 } - \varpi } z.
\]
В результате получим векторное дифференциальное уравнение вида
\begin{gather}
\dfrac{ dz }{ dt } = P(t) z + q(t) + x( t, z )
\ \ \left( \left( t, z \right) \in D_3 \right),
                                                    \label{DEV:solution->0}
                                                                        \\
P(t) := \Re \, D(t) + \left( \varpi - \frac{ \rho_{ k+1 } }{ \varrho_1 } \right)
\left( \varepsilon_1(t) \right)^{ -1 } \dfrac{ d \varepsilon_1(t) }{ dt } E,
\ \ q(t) := \left( \varepsilon_1(t) \right)^
{ \varpi - \frac { \rho_{ k+1 } } { \varrho_1 } }
Q^{ -1 }(t) g(t),
                                                            \notag
                                                                        \\
x( t, z ) := \left( \varepsilon_1(t) \right)^
{ \varpi - \frac { \rho_{ k+1 } } { \varrho_1 } }
Q^{ -1 }(t) \sum_{ l=1 }^m \varepsilon_l(t) \psi_l \left( t,
\left( \varepsilon_1(t) \right)^{ \frac { \rho_{ k+1 } } { \varrho_1 } -
\varpi } Q(t) z \right),
                                                            \label{x(t,z)}
\end{gather}
область
$ D_3 := \left\{ ( t, z ) \ \left| \ t \in I,
\ z \in \mathbb R^{ n \times 1 }, \ || z || \le \tilde a \right. \right\}, $
число
\[
\tilde a := \hat a_0 \left( \sup\limits_{ t \in I } \varepsilon_1(t) \right)^{
\varpi - \frac { \rho_{ k+1 } } { \varrho_1 } }
\left( n \sup\limits_{ t \in I } || Q(t) || \right)^{ -1 }.
\]

Проверим, выполнено ли для вектор-функции $ x( t, z ) $ в области
$ D_3 $ по переменной $ z $ условие Липшица.
Для оценки сверху нормы $ || x( t, z_1 ) - x( t, z_2 ) || $
$ \left( \forall \ ( t, z_1 ), ( t, z_2 ) \in D_3 \right) $ воспользуемся
аналогом формулы Лагранжа конечных приращений скалярных функций.
Таким образом, достаточно оценить сверху норму матрицы Якоби вектор--функции
$ x( t, z ) $ по переменной $ z $.
Найдем эту матрицу, используя формулу \eqref{x(t,z)}:
$$
\frac{ \partial x( t, z ) }{ \partial z } = Q^{ -1 }(t)
\sum_{ l=1 }^m \varepsilon_l(t) \frac{ \partial \psi_l \left( t,
\left( \varepsilon_1(t) \right)^{ \frac { \rho_{ k+1 } } { \varrho_1 } -
\varpi } Q(t) z \right) }{ \partial r } Q(t).
$$
Следовательно,
$$
\left\| \frac{ \partial x( t, z ) }{ \partial z } \right\| =
O( \varepsilon_1(t) ) \sum_{ l=1 }^m \left\|
\frac{ \partial \psi_l \left( t, \left( \varepsilon_1(t)
\right)^{ \frac{ \rho_{ k+1 } }{ \varrho_1 } - \varpi } Q(t) z
\right) }{ \partial r } \right\|
\ \ \left( t \to + \infty \right).
$$
В результате подстановки свойств \eqref{Est:JM_psi_l(t,r)} в это неравенство,
получим искомую оценку матрицы Якоби
$$
\left\| \frac{ \partial x( t, z ) }{ \partial z } \right\| =
O \left( \left( \varepsilon_1(t)
\right)^{ 1 + \frac{ \rho_{ k+1 } }{ \varrho_1 } - \varpi } \right)
\ \ \left( t \to + \infty, \ || z || \le \tilde a \right).
$$
В итоге получили требуемое условие Липшица для вектор-функции
$ x( t, z ) $ в области $ D_3 $ по переменной $ z $:
\begin{multline}
|| x( t, z_1 ) - x( t, z_2 ) || \le
x^*(t) || z_1 - z_2 ||,
\ \ 0 < x^*(t) := const \left( \varepsilon_1(t) \right)^
{ 1 + \frac { \rho_{ k+1 } } { \varrho_1 } - \varpi }
                                                                        \\
\left( \forall \ ( t, z_1 ), ( t, z_2 ) \in D_3,
\ t \to + \infty \right).
                                                            \label{IE:Lip}
\end{multline}
Так как выполнено неравенство
$ \rho_{ k+1 } \ge \varrho_1 ( 2 \varpi - 1 ) $, то функция
$ x^*(t) \in \mathcal R_1(I) $.

В силу тождеств \eqref{Id:psi_l(t,0)=0} и равенства \eqref{x(t,z)}, очевидно,
что вектор--функция $ x( t, 0 ) \equiv 0 $ $ ( t \in I ) $.
Поэтому, учитывая условие Липшица \eqref{IE:Lip}, заметим, что для
нелинейности $ x(t, z) $ в области $ D_3 $ будет иметь место оценка
\begin{equation}
|| x(t, z) || \le x^*(t) || z || \ \ \left( t \to + \infty \right).
                                                        \label{IE:estimateNL}
\end{equation}

Запишем векторное дифференциальное уравнение \eqref{DEV:solution->0} в
скалярной форме
\begin{gather}
\dfrac{ d z_j }{ dt } = q_j(t) + p_{ jj }(t) z_j + x_j(t, z_1, \dots, z_n)
\ \ \left( j = \overline{ 1, n } \right),
                                                    \label{SDE:solution->0}
                                                                        \\
z_j := ( z )_j, \ q_j(t) := ( q(t) )_j,
\ p_{ jj }(t) := ( P(t) )_{ jj },
\ x_j(t, z_1, \dots, z_n) := ( x( t, z ) )_j.
                                                            \label{p_jj(t)}
\end{gather}
Применим для нахождения ограниченных при $ t \ge t_0 $ решений системы
дифференциальных уравнений \eqref{SDE:solution->0} специальный метод
последовательных приближений, аналогичный тому, который использован в ходе
доказательства теоремы 1.1 из \S \ 1 главы IV кандидатской диссертации
\cite{Bib:Kostin_PHDT} (страница 67).
Пусть $ z_{ j \, s-1 }(t) $ ($ j = \overline{ 1, n } $) означает
$ ( s-1 ) $-е приближение, а $ z_{ js }(t) $ --- $ s $-е
($ j = \overline{ 1, n } $).
Положим $ z_{ 10 }(t) := \dots := z_{ n0 }(t) := 0. $
Определим $ s $-е приближение из системы
\begin{equation}
\dfrac{ d z_{ js } }{ dt } = q_j(t) + x_j(t, z_{ 1 \, s-1 }, \dots,
z_{ n \, s-1 } ) + p_{ jj }(t) z_{ js }
\ \ \left( j = \overline{ 1, n } \right),
                                            \label{SDE:Sth-approximation}
\end{equation}
выбирая начальные значения для функций $ z_{ js }(t) $
($ j = \overline{ 1, n } $) так, чтобы эти функции выражались через
$ z_{ j \, s-1 }(t) $ ($ j = \overline{ 1, n } $) формулами вида
\begin{multline}
z_{ js }(t) = \int\limits_{ a_j }^t q_j( \tau )
\exp \int\limits_\tau^t p_{ jj }(t) \, dt \, d \tau +
z_j(t_0) \exp \int\limits_{ t_0 }^t p_{ jj }( \tau ) \, d \tau +
                                                                        \\
+ \int\limits_{ a_j }^t x_j( \tau, z_{ 1 \, s-1 }( \tau ), \dots,
z_{ n \, s-1 } ( \tau ) ) \exp \int\limits_\tau^t p_{ jj }(t) \, dt \, d \tau
\ \ ( j = \overline{ 1, n} ),
                                                        \label{solution->0}
\end{multline}
где каждый предел интегрирования $ a_j $ равен либо $ t_0 $, либо
$ + \infty $; начальные значения $ z_j(t_0) := 0 $ для тех индексов
$ j $, для которых не выполняются условия \eqref{No_of_parameters_O}.
Нетрудно проверить простым дифференцированием, что при любом выборе указанных
пределов интегрирования, равенства \eqref{solution->0} будут давать нам
некоторое частное решение системы дифференциальных уравнений
\eqref{SDE:Sth-approximation} (если только величины, входящие в
\eqref{solution->0} не теряют смысла).
Пределы интегрирования $ a_j $ будем выбирать так:
\[
a_j :=
\begin{cases}
+ \infty, & \text{если }
p_{ jj }(t) - \Re \, \varsigma_j(t) \ge 0 \ (t \ge t_0),                 \\
t_0, & \text{если }
p_{ jj }(t) - \Re \, \varsigma_j(t) \le 0 \ (t \ge t_0)
\end{cases}
\ \ ( j = \overline{ 1, n } ).
\]

Постараемся сделать так, чтобы все последовательные
приближения были ограничены по модулю одним и тем же числом.
С этой целью предположим, что
\begin{equation}
| z_{ j \, s-1 }(t) | \le \epsilon_0
\ \ ( j = \overline{ 1, n } ),
                                                        \label{IE:constraint}
\end{equation}
где $ \epsilon_0 $ --- некоторая константа,
$ 0 < \epsilon_0 \le \tilde a $ и потребуем, чтобы такие же
неравенства имели место и для $ s $-го приближения.
Принимая во внимание свойство \eqref{IE:estimateNL},
нетрудно заметить, что при выполнении неравенств \eqref{IE:constraint},
функции $ z_{ js }(t) $ ($ j = \overline{ 1, n } $) будут мажорироваться по
модулю функциями $ \xi_j(t, t_0, \epsilon_0) $
($ j = \overline{ 1, n } $), которые определяются из
равенств
\begin{multline}
\xi_j(t, t_0, \epsilon_0) :=
b_j \int\limits_{ a_j }^t | q_j( \tau ) |
\exp \int\limits_{ \tau }^t p_{ jj }(t) \, dt \, d \tau + z_j^*(t) +
																																				\\
+ \epsilon_0 b_j \int\limits_{ a_j }^t x^*( \tau )
\exp \int\limits_{ \tau }^t p_{ jj }(t) \, dt \, d \tau
\ \ ( j = \overline{ 1, n } ),
                                            \label{xi_j(t t_0 epsilon_0)}
\end{multline}
где функции
\[
z_j^*(t) := | z_j(t_0) | \exp \int\limits_{ t_0 }^t p_{ jj }( \tau )
\, d \tau \ \ \left( j = \overline{ 1, n } \right),
\]
числа
\[
b_j :=
\begin{cases}
\phantom{ - } 1,      & \text{если } a_j = t_0,                         \\
-1,     & \text{если } a_j = + \infty
\end{cases}
\ \ ( j = \overline{ 1, n } ).
\]
Неравенства \eqref{IE:constraint} заведомо будут выполняться
для $ s $-го приближения, если эти неравенства имеют место
для функций $ \xi_j(t, t_0, \epsilon_0) $
($ j = \overline{ 1, n } $).
Запишем неравенства
\[
\xi_j(t, t_0, \epsilon_0) \le \epsilon_0
\ \ ( j = \overline{ 1, n } )
\]
и постараемся найти с их помощью константу $ \epsilon_0 $.
Замечая, что функции $ \xi_j(t, t_0, \epsilon_0) $
($ j = \overline{ 1, n } $) линейны относительно постоянной
$ \epsilon_0 $, представим последние неравенства в виде
\[
a_j(t, t_0) + \epsilon_0 b_j(t, t_0) \le \epsilon_0
\ \ ( j = \overline{ 1, n } ),
\]
где
\begin{align*}
& a_j(t, t_0) := b_j \int\limits_{ a_j }^t | q_j( \tau ) |
\exp \int\limits_{ \tau }^t p_{ jj }(t) \, dt \, d \tau + z_j^*(t),
                                                                        \\
& b_j(t, t_0) := b_j \int\limits_{ a_j }^t x^*( \tau )
\exp \int\limits_{ \tau }^t p_{ jj }(t) \, dt \, d \tau
\ \ \left( j = \overline{ 1, n } \right).
\end{align*}
Искомое $ \epsilon_0 $ легко находится из последних неравенств, если
$ b_j(t, t_0) < 1 $ при $ t \ge t_0 $
$ \left( j = \overline{ 1, n } \right) $ и если, кроме того,
\[
\frac{ a_j(t, t_0) }{ 1 - b_j(t, t_0) } \le \tilde a, \ \ t \ge t_0
\ \ \left( j = \overline{ 1, n } \right),
\]
причём в этом случае в качестве $ \epsilon_0 $ можно взять константу
\[
\epsilon_0 := \max_{ j = \overline{ 1, n } } \left(
\sup_{ t \in I } \frac{ a_j(t, t_0) }{ 1 - b_j(t, t_0) } \right).
\]

В силу условий теоремы функции $ q_j(t) $,
$ x^*(t) \in \mathcal R_1(I) $ $ \left( j = \overline{ 1, n } \right) $.
Значит, для выражений $ a_j(t, t_0) - z_j^*(t) $, $ b_j(t, t_0) $
$ \left( j = \overline{ 1, n } \right) $ выполнены требования лемм 3 и 3$'$
(с учётом замечаний к ним) из \S \, 2 главы II кандидатской диссертации
\cite{Bib:Kostin_PHDT} (страницы 46 и 47 соответственно).
Следовательно, выражения $ a_j(t, t_0) $, $ b_j(t, t_0) $
$ \left( j = \overline{ 1, n } \right) $ ограничены на промежутке $ I $.
Так как подынтегральные функции во внешних интегралах в выражениях
$ a_j(t, t_0) - z_j^*(t) $, $ b_j(t, t_0) $
$ \left( j = \overline{ 1, n } \right) $ неотрицательны, то за счёт выбора
числа $ t_0 $ достаточно большим, а $ | z_j( t_0 ) | $  (индексы $ j $ такие,
что выполнены условия \eqref{No_of_parameters_O}, $ z_j(t_0) \in \mathbb R $)
достаточно малыми можем сделать величины
$ \sup\limits_{ t \in I } a_j(t, t_0) $ и
$ \sup\limits_{ t \in I } b_j(t, t_0) $
$ \left( j = \overline{ 1, n } \right) $ сколь угодно малыми.
Следовательно, для достаточно большого числа $ t_0 $ и, либо малых
$ | z_j( t_0 ) | $ (индексы $ j $ такие, что выполнены условия
\eqref{No_of_parameters_O}), либо начальных значений $ z_j( t_0 ) \equiv 0 $
(в случае отсутствия таковых индексов $ j $), величины
\begin{equation*}
\max\limits_{ j = \overline{ 1, n } } \left( \sup\limits_{ t \in I }
b_j(t, t_0) \right) < 1,
\ \ \epsilon_0 \le \tilde a
\end{equation*}
(где $ \tilde a $ определяет область $ D_3 $).
Таким образом, для системы \eqref{xi_j(t t_0 epsilon_0)} выполнены все
условия теоремы 1.1 из \S \ 1 главы IV кандидатской диссертации
\cite{Bib:Kostin_PHDT} (страница 67).
Поэтому для достаточно большого числа $ t_0 $ система дифференциальных
уравнений \eqref{SDE:solution->0} будет заведомо иметь хотя бы одно
вещественное ограниченное при $ t \ge t_0 $ частное решение $ z_j(t) $
$ \left( j = \overline{ 1, n } \right) $ с условием
$ | z_j(t) | \le \epsilon_0 $
$ \left( t \in I; \ j = \overline{ 1, n } \right) $.
Более того, решение $ z_j(t) $ $ \left( j = \overline{ 1, n } \right) $
зависит от стольких достаточно малых по абсолютной величине произвольных
постоянных, сколько имеется индексов $ j $ таких, что выполнены условия
\eqref{No_of_parameters_O}.
Возвращаясь обратно к вектор--функции $ r $, получим требуемое.

Как известно, общее решение системы обыкновенных дифференциальных уравнений
$ n $--го порядка зависит от $ n $ произвольных параметров, следовательно в
семействе решений вида \eqref{CV:r_j} могут фигурировать не более чем $ n $
независимых скалярных постоянных.

\end{proof}

В следующем замечании указаны выражения для коэффициентов
$ \gamma_{ j \wp } $ в первых слагаемых сумм в функциях $ d_j(t) $, которые
могут понадобиться при применении теорем (\ref{T:AsymptChart_O} --
\ref{T:AsymptChart}) об асимптотическом характере с целью вычисления
количества независимых скалярных параметров, которые фигурируют в погрешности
$ r $ из семейства решений вида \eqref{CV:r_j} основного векторного
дифференциального уравнения \eqref{DEV:main}.

\begin{remark}

Рассматривая ход доказательства теоремы 1 из работы авторов \cite{Bib:Odessa01},
заметим, что числа
$ \gamma_{ j \wp } = \operatorname{ M }
\left( \left( P_0^{ -1 } A_\wp(t) P_0 \right)_{ jj } \right) $,
где матрица $ P_0 \in \mathbb C^{ n \times n } $ составлена из собственных
векторов матрицы $ A $, $ j = \overline{ 1, n } $, $ | \wp | = 1 $,
$ \wp \in \mathbb N_0^m $.
Подставляя в эти равенства тождества \eqref{A_p(t), |p|=1}, получим искомые
коэффициенты:
\begin{equation}
\gamma_{ j \wp } = \operatorname{ M } \left( \left( P_0^{ -1 }
\frac { \partial } { \partial { y } } f_l ( t, \varphi_0( t, c_0 ) )
P_0 \right)_{ jj } \right),
\ j = \overline{ 1, n },
\ \wp = ( \underbrace{ 0, \dots, 0 }_{ l-1 }, 1, 0, \dots, 0 ),
                                                \label{gamma_jp, |p|=1}
\end{equation}
где мультииндексы $ \wp \in \mathbb N_0^m $, функции
$ \varepsilon_l(t) \notin \mathcal R_1(I) $.

\end{remark}

В следующей теореме исследован асимптотический характер формального частного
решения (семейства решений) \eqref{FSolution} основного вещественного
векторного квазилинейного обыкновенного дифференциального уравнения
\eqref{DEV:main} в том случае, когда матрица $ A $ не имеет кратных
собственных значений, а оценка погрешности $ r $ получена в виде $ o $.

\begin{theorem}
                                                    \label{T:AsymptChart_o}
Пусть выполнены все условия теоремы \ref{T:AsymptChart_O} и либо функции
\begin{equation}
p_{ jj }(t) - \Re \, \varsigma_j(t) \ge 0 \ \ ( t \in I ),
                                                \label{C:No_parameters_o}
\end{equation}
либо справедливы требования
\begin{equation}
\int\limits_I \left( p_{ jj }(t) - \Re \, \varsigma_j(t) \right) \, dt =
- \infty
                                                \label{No_of_parameters_o}
\end{equation}
$ \left( j = \overline{ 1, n } \right). $
Тогда для достаточно большого числа $ t_0 $ у векторного дифференциального
уравнения \eqref{DEV:main} на промежутке $ I $ существует хотя бы одно
частное решение вида \eqref{CV:r_j}, где вектор--функции
$ \varphi_0( t, c_0 ), \varphi_{ sp }( t, c_{ sp } ) \in
\mathcal K_3^{ n \times 1 }(A) \cap C^1( I ) $
($ s = \overline{ 1, k } $, $ p = \overline{ 1, \varkappa_s } $), причём
погрешность
\begin{equation*}
r = o \left( \left( \varepsilon_1(t) \right)^{
\frac { \rho_{ k+1 } } { \varrho_1 } - \varpi } \right)
\ \ ( t \to + \infty ).
\end{equation*}
Более того, на промежутке $ [ t_1, + \infty ) $ погрешность $ r $ зависит от
стольких произвольных скалярных постоянных, сколько имеется индексов
$ j \in \left\{ 1, 2, 3, \dots, n \right\} $ таких, что выполнены условия
\eqref{No_of_parameters_o} ($ t_1 \in \mathbb R $ --- достаточно большое
число, определяется этими постоянными, $ t_1 \ge t_0 $).
(В семействе решений вида \eqref{CV:r_j} могут фигурировать не более чем
$ n $ независимых скалярных параметров.)

\end{theorem}

\begin{proof}
В ходе доказательства теоремы \ref{T:AsymptChart_O} было показано, что при
условиях теоремы для достаточно большого числа $ t_0 $ система
дифференциальных уравнений \eqref{SDE:solution->0} будет заведомо иметь хотя
бы одно вещественное ограниченное при $ t \ge t_0 $ частное решение
$ z_j(t) $ $ \left( j = \overline{ 1, n } \right) $ с условием
$ | z_j(t) | \le \epsilon_0 $
$ \left( t \in I; \ j = \overline{ 1, n } \right) $.
Там же было показано, что при условиях теоремы функции $ q_j(t) $,
$ x^*(t) \in \mathcal R_1(I) $ $ \left( j = \overline{ 1, n } \right) $.
Поэтому, учитывая, что справедливы либо неравенства
\eqref{C:No_parameters_o}, либо предположения \eqref{No_of_parameters_o},
заметим, что для выражений $ a_j(t, t_0) - z_j^*(t) $, $ b_j(t, t_0) $
$ \left( j = \overline{ 1, n } \right) $ выполнены требования лемм 3 и 3$'$
(с учётом замечаний к ним) из \S \, 2 главы II кандидатской диссертации
\cite{Bib:Kostin_PHDT} (страницы 46 и 47 соответственно).
Следовательно, выражения $ a_j(t, t_0) - z_j^*(t) $, $ b_j(t, t_0) \to 0 $
$ \left( t \to + \infty; \ j = \overline{ 1, n } \right) $.

Пусть теперь в последовательных приближениях \eqref{solution->0} начальные
значения $ z_j(t_0) := 0 $ для тех индексов
$ j \in \left\{ 1, 2, 3, \dots, n \right\} $, для которых выполняются
неравенства \eqref{C:No_parameters_o}.
Очевидно, что в таком случае для любых начальных значений $ z_j(t_0) $
(индексы $ j \in \left\{ 1, 2, 3, \dots, n \right\} $ такие, что справедливы
требования \eqref{No_of_parameters_o}) пределы
$ \lim\limits_{ t \to + \infty } z_j^*(t) = 0 $
$ \left( j = \overline{ 1, n } \right). $
Значит, найдется такое достаточно большое число $ t_1 \in \mathbb R $
($ t_1 \ge t_0 $), что величины $ \sup\limits_{ t \ge t_1 } a_j(t, t_0) $ и
$ \sup\limits_{ t \ge t_1 } b_j(t, t_0) $
$ \left( j = \overline{ 1, n } \right) $ будут сколь угодно малыми.
Таким образом, для системы \eqref{xi_j(t t_0 epsilon_0)} на промежутке
$ [ t_1, + \infty ) $ выполнены все условия теоремы 1.1 из \S \ 1 главы IV
кандидатской диссертации \cite{Bib:Kostin_PHDT} (страница 67).
Поэтому для любых начальных значений $ z_j(t_0) \in \mathbb R $ (индексы
$ j \in \left\{ 1, 2, 3, \dots, n \right\} $ такие, что справедливы
требования \eqref{No_of_parameters_o}) на промежутке $ [ t_1, + \infty ) $
решение $ z_j(t) $ $ \left( j = \overline{ 1, n } \right) $ зависит от
стольких произвольных постоянных, сколько имеется индексов
$ j \in \left\{ 1, 2, 3, \dots, n \right\} $ таких, что выполнены условия
\eqref{No_of_parameters_o}.

Так как функции $ z_{ js }(t) $
$ \left( j = \overline{ 1, n }, \ s \in \mathbb N_0 \right) $ на промежутке
$ I $ мажорируются по модулю функциями $ \xi_j(t, t_0, \epsilon_0) $
($ j = \overline{ 1, n } $), которые стремятся к нулю при $ t \to + \infty $,
то $ z_{ js }(t) = o(1) $
$ \left( j = \overline{ 1, n }, \ s \in \mathbb N_0,
\ t \to + \infty \right) $.
В итоге для системы \eqref{SDE:solution->0} и равенств \eqref{solution->0}
на промежутке $ [ t_1, + \infty ) $ выполнены все условия теоремы 1.2 из \S
\ 1 главы IV кандидатской диссертации \cite{Bib:Kostin_PHDT} (страница 68).
Поэтому для достаточно большого числа $ t_0 $ система дифференциальных
уравнений \eqref{SDE:solution->0} будет заведомо иметь хотя бы одно
вещественное частное решение $ z_j(t) $
$ \left( j = \overline{ 1, n } \right) $ с условием $ z_j(t) = o(1) $
$ \left( j = \overline{ 1, n }; \ t \to + \infty \right) $.
Более того, на промежутке $ [ t_1, + \infty ) $ решение $ z_j(t) $
$ \left( j = \overline{ 1, n } \right) $ зависит от стольких произвольных
постоянных, сколько имеется индексов
$ j \in \left\{ 1, 2, 3, \dots, n \right\} $ таких, что выполнены условия
\eqref{No_of_parameters_o} ($ t_1 \in \mathbb R $ --- достаточно большое
число, определяется этими постоянными, $ t_1 \ge t_0 $).
Возвращаясь обратно к вектор--функции $ r $, получим требуемое.

Как известно, общее решение системы обыкновенных дифференциальных уравнений
$ n $--го порядка зависит от $ n $ произвольных параметров, следовательно в
семействе решений вида \eqref{CV:r_j} могут фигурировать не более чем $ n $
независимых скалярных постоянных.

\end{proof}

В следующей теореме получен более высокий порядок малости в оценке
погрешности $ r $ по сравнению с теоремой \ref{T:AsymptChart_O}.

\begin{theorem}

Пусть выполнены все условия теоремы \ref{T:AsymptChart_O} для числа
$ k $ такого, что
$ \rho_{ k+1 } \ge \rho_{ k_0 + 1 } + \varpi \varrho_1 $
($ k_0 \in \mathbb N $).
Тогда для достаточно большого числа $ t_0 $ у векторного дифференциального
уравнения \eqref{DEV:main} на промежутке $ I $ существует хотя бы одно
частное решение вида
                                            \label{T:AsymptChart}
\begin{equation}
y = \varphi_0( t, c_0 ) + \sum_{ s=1 }^{ k_0 }
\sum_{ p=1 }^{ \varkappa_s } \nu_{ sp }(t)
\varphi_{ sp }( t, c_{ sp } ) + r,
\ \ r = O \left( \left( \varepsilon_1(t) \right)^{
\frac { \rho_{ k_0 + 1 } } { \varrho_1 } } \right)
\ ( t \to + \infty ),
                                            \label{AsympAproxSolution}
\end{equation}
где вектор--функции
$ \varphi_0( t, c_0 ), \varphi_{ sp }( t, c_{ sp } ) \in
\mathcal K_3^{ n \times 1 }(A) \cap C^1( I ) $
($ s = \overline{ 1, k_0 } $, $ p = \overline{ 1, \varkappa_s } $).
Более того, погрешность $ r $ зависит от стольких достаточно малых по
абсолютной величине произвольных скалярных постоянных, сколько имеется
индексов $ j $ таких, что справедливы условия \eqref{No_of_parameters_O}.
(В семействе решений \eqref{AsympAproxSolution} могут фигурировать не более
чем $ n $ независимых скалярных параметров.)

\end{theorem}

\begin{proof}
В силу теоремы \ref{T:AsymptChart_O} у векторного дифференциального
уравнения \eqref{DEV:main} существует хотя бы одно частное решение вида
\eqref{CV:r_j}, причём для погрешности $ r $ будет справедлива оценка
\eqref{r_j = O}.
Далее, из условий теоремы следует, что будет выполнено
следующее свойство
\[
r = O \left( \left( \varepsilon_1(t)
\right)^{ \frac { \rho_{ k_0 + 1 } } { \varrho_1 } } \right)
\ \ ( t \to + \infty ).
\]
Разобьём векторную сумму $ s( t ) $ на две части.
Пусть в первой будут слагаемые содержащие функции $ \nu_{ sp }(t) $ ранга
не превосходящего $ \rho_{ k_0 } $, а вторую оценим с помощью свойств
\eqref{estimate_nu}
\[
\sum_{ s = k_0 + 1 }^k \sum_{ p=1 }^{ \varkappa_s } \nu_{ sp }(t)
\varphi_{ sp }( t, c_{ sp } ) =
O \left( \left( \varepsilon_1(t)
\right)^{ \frac { \rho_{ k_0 + 1 } } { \varrho_1 } } \right)
\ \ ( t \to + \infty ).
\]
Следовательно, представили частное решение векторного дифференциального
уравнения \eqref{DEV:main} в требуемом виде.

\end{proof}

\section{Примеры}

В завершении статьи приведем три простых примера.
Следующий пример иллюстрирует применение теоремы \ref{T:NPI_AsymptChart_O} об
асимптотическом характере формального частного решения (семейства решений)
\eqref{FSolution} основного вещественного векторного квазилинейного
обыкновенного дифференциального уравнения \eqref{DEV:main} в том случае,
когда матрица $ A $ не имеет чисто мнимых собственных значений, а класс
$ \mathcal K_1(A) $ состоит из периодических функций.

\begin{example}

Пусть числа $ t_0 > 1 $, $ n := m := 2 $, $ k := 1 $,
$ \tau \in \mathbb R^+ $ фиксировано, матрица
$ A \in \mathbb R^{ 2 \times 2 } $, $ \Re \, \lambda_j(A) \ne 0 $
$ \left( j = \overline{ 1, 2 } \right) $.
Зададим класс $ \mathcal K_1(A) $ так, как указано в примере
\ref{Ex:PeriodicK1}, а ранги степенно-логарифмических функций как в примере
\ref{E:rang}.

Рассмотрим в области $ D := \left\{ ( t, { y } ) \ \left| \ t \in I,
\ { y } \in \mathbb R^{ 2 \times 1 }, \right.
|| { y } - \varphi_0( t ) || \le a \right\} $
(число $ a \in \mathbb R^+ $) векторное дифференциальное уравнение
\begin{equation}
\frac{ dy }{ dt } = A y + f(t) + \frac{1}{t}
\begin{pmatrix}
  ( y )_1 + ( y )_2^2 g_1(t) \\
  ( y )_2 \\
\end{pmatrix}
+ \frac{1}{ t \ln t }
\begin{pmatrix}
  ( y )_2 \\
  ( y )_1 + ( y )_1 ( y )_2 g_2(t) \\
\end{pmatrix},
                                                    \label{DEV:NPI_example}
\end{equation}
причём вектор--функция
$ \varphi_0( t ) : \mathbb R \to \mathbb R^{ 2 \times 1 } $ ---
единственное $ \tau $\nobreakdash-\hspace{0mm}периодическое частное решение
из класса $ C^1( \mathbb R ) $ укороченного векторного дифференциального
уравнения соответствующего уравнению \eqref{DEV:NPI_example}:
\begin{equation*}
\frac{ d \varphi_0( t ) }{ dt } = A \varphi_0( t ) + f(t),
\end{equation*}
вектор--функция $ f( t ) : \mathbb R \to \mathbb R^{ 2 \times 1 } $ и функции
$ g_1( t ), \ g_2( t ) : \mathbb R \to \mathbb R $ --- некоторые наперед
заданные $ \tau $\nobreakdash-\hspace{0mm}периодические из класса
$ C^1( \mathbb R ) $.

Очевидно, что в рассматриваемом случае
$$
\rho_s = s \ ( s \in \mathbb N ), \ \varkappa_1 = 2,
\ \nu_{ 11 }(t) \equiv \varepsilon_1(t) \equiv \frac{1}{t},
\ \nu_{ 12 }(t) \equiv \varepsilon_2(t) \equiv \frac{1}{ t \ln t }.
$$

Для нахождения вектор--функций
$ \varphi_{ 11 }(t), \ \varphi_{ 12 }(t) :
\mathbb R \to \mathbb R^{ 2 \times 1 } $
воспользуемся замечанием \ref{R:first_terms}.
Найдём их как единственные $ \tau $\nobreakdash-\hspace{0mm}периодические
частные решения из класса $ C^1( \mathbb R ) $ векторных дифференциальных
уравнений из совокупности
\begin{equation*}
\left[
\begin{array}{l}
\dfrac{ d \varphi_{ 11 }(t) }{ dt } = A \varphi_{ 11 }(t) +
\begin{pmatrix}
  ( \varphi_0( t ) )_1 + \left( \varphi_0( t ) \right)_2^2 g_1(t) \\
  ( \varphi_0( t ) )_2 \\
\end{pmatrix},
                                                                        \\
\dfrac{ d \varphi_{ 12 }(t) }{ dt } = A \varphi_{ 12 }(t) +
\begin{pmatrix}
  ( \varphi_0( t ) )_2 \\
  ( \varphi_0( t ) )_1 + \left( \varphi_0( t ) \right)_1
\left( \varphi_0( t ) \right)_2 g_2(t) \\
\end{pmatrix}.
\end{array}
\right.
\end{equation*}

И так, выполнены все условия теоремы \ref{T:NPI_AsymptChart_O},
следовательно, для достаточно большого числа $ t_0 $ у векторного
дифференциального уравнения \eqref{DEV:NPI_example} на промежутке $ I $
существует хотя бы одно частное решение вида
$$
y = \varphi_0( t ) + \frac{1}{t} { \varphi }_{ 11 }(t) +
\frac{1}{ t \ln t } { \varphi }_{ 12 }(t) + r,
\ \ r = O \left( \frac{ 1 }{ t^2} \right) \ ( t \to + \infty ).
$$
Более того, на промежутке $ [ t_1, + \infty ) $ погрешность $ r $ зависит от
стольких произвольных скалярных постоянных, сколько имеется индексов $ j $
таких, что справедливы условия
$ j \in \left\{ 1, 2 \right\}, \ \Re \, \lambda_j(A) < 0, $
$ t_1 \in \mathbb R $ --- достаточно большое число, определяется этими
постоянными ($ t_1 \ge t_0 $).

\end{example}

Следующий пример иллюстрирует применение теоремы \ref{T:NPI_AsymptChart_O} об
асимптотическом характере формального частного решения (семейства решений)
\eqref{FSolution} основного вещественного векторного квазилинейного
обыкновенного дифференциального уравнения \eqref{DEV:main} в том случае,
когда матрица $ A $ не имеет чисто мнимых собственных значений, а класс
$ \mathcal K_1(A) $ задан так, как указано в примере
\ref{Ex:Decreasing_Exp_LC_K1}.

\begin{example}

Пусть числа $ t_0 > 1 $, $ n := m := 2 $, $ k := 1 $, матрица
$ A \in \mathbb R^{ 2 \times 2 } $, $ \Re \, \lambda_j(A) \ne 0 $
$ \left( j = \overline{ 1, 2 } \right) $ и выполнено условие
\eqref{C:exponent_EV_separated}, где множество
\begin{equation*}
\Gamma_3 := \left\{ \left. \sum_{ s=1 }^4 k_s \omega_s \, \right|
k_s \in \mathbb N_0, \ \omega_s \in \Omega_2
\ \left( s = \overline{ 1, 4 } \right) \right\},
\ \ \Omega_2 := \left\{ -1, \ -\sqrt 2, \ -e, \ -\pi \right\}.
\end{equation*}
В качестве класса $ \mathcal K_1(A) $ возьмём множество конечных линейных
комбинаций \eqref{Set:Decreasing_Exp_LC}.
Зададим ранги степенно-логарифмических функций так, как указано в примере
\ref{E:rang}.

Рассмотрим в области $ D := \left\{ ( t, { y } ) \ \left| \ t \in I,
\ { y } \in \mathbb R^{ 2 \times 1 }, \right.
|| { y } - \varphi_0( t ) || \le a \right\} $
(число $ a \in \mathbb R^+ $) векторное дифференциальное уравнение
\begin{multline}
\frac{ dy }{ dt } = A y + m_0 +
\begin{pmatrix}
  e^{ -t } \\
  e^{ -e t } \\
\end{pmatrix}
+ \frac{1}{t} \left( m_1 +
\begin{pmatrix}
  ( y )_1 + ( y )_2^2 e^{ - \pi t } \\
  ( y )_2 \\
\end{pmatrix}
\right) +
                                                                        \\
+ \frac{1}{ t \ln t } \left( m_2 +
\begin{pmatrix}
  ( y )_1 \\
  ( y )_2 + ( y )_1 ( y )_2 e^{ - \sqrt{2} \, t } \\
\end{pmatrix}
\right),
                                                \label{DEV:NPI_DELC_example}
\end{multline}
причём вектор--функция $ \varphi_0( t ) $ --- единственное частное решение из
класса $ \mathcal K_1^{ 2 \times 1 }(A) \cap C^1( I ) $ укороченного
векторного дифференциального уравнения соответствующего уравнению
\eqref{DEV:NPI_DELC_example}:
\begin{equation*}
\frac{ d \varphi_0( t ) }{ dt } = A \varphi_0( t ) + m_0 +
\begin{pmatrix}
  e^{ -t } \\
  e^{ -e t } \\
\end{pmatrix},
\end{equation*}
векторы $ m_0, \ m_1, \ m_2 \in \mathbb R^{ 2 \times 1 } $.

Очевидно, что в рассматриваемом случае
$$
\rho_s = s \ ( s \in \mathbb N ), \ \varkappa_1 = 2,
\ \nu_{ 11 }(t) \equiv \varepsilon_1(t) \equiv \frac{1}{t},
\ \nu_{ 12 }(t) \equiv \varepsilon_2(t) \equiv \frac{1}{ t \ln t }.
$$

Для нахождения вектор--функций $  \varphi_{ 11 }(t) $,
$  \varphi_{ 12 }(t) $ воспользуемся замечанием \ref{R:first_terms}.
Найдём их как единственные частные решения из класса
$ \mathcal K_1^{ 2 \times 1 }(A) \cap C^1( I ) $ векторных
дифференциальных уравнений из совокупности
\begin{equation*}
\left[
\begin{array}{l}
\dfrac{ d \varphi_{ 11 }(t) }{ dt } = A \varphi_{ 11 }(t) + m_1 +
\begin{pmatrix}
  ( \varphi_0( t ) )_1 + \left( \varphi_0( t ) \right)_2^2 e^{ - \pi t } \\
  ( \varphi_0( t ) )_2 \\
\end{pmatrix},
                                                                        \\
\dfrac{ d \varphi_{ 12 }(t) }{ dt } = A \varphi_{ 12 }(t) + m_2 +
\begin{pmatrix}
  ( \varphi_0( t ) )_1 \\
  ( \varphi_0( t ) )_2 + \left( \varphi_0( t ) \right)_1
\left( \varphi_0( t ) \right)_2 e^{ - \sqrt{2} \, t } \\
\end{pmatrix}.
\end{array}
\right.
\end{equation*}

И так, выполнены все условия теоремы \ref{T:NPI_AsymptChart_O},
следовательно, для достаточно большого числа $ t_0 $ у векторного
дифференциального уравнения \eqref{DEV:NPI_DELC_example} на промежутке $ I $
существует хотя бы одно частное решение вида
\begin{equation}
y = \varphi_0( t ) + \frac{1}{t} { \varphi }_{ 11 }(t) +
\frac{1}{ t \ln t } { \varphi }_{ 12 }(t) + r,
\ \ r = O \left( \frac{ 1 }{ t^2} \right) \ ( t \to + \infty ).
                                            \label{NPI_DELC_exampleSolution}
\end{equation}
Более того, на промежутке $ [ t_1, + \infty ) $ погрешность $ r $ зависит от
стольких произвольных скалярных постоянных, сколько имеется индексов $ j $
таких, что справедливы условия
$ j \in \left\{ 1, 2 \right\}, \ \Re \, \lambda_j(A) < 0, $
$ t_1 \in \mathbb R $ --- достаточно большое число, определяется этими
постоянными ($ t_1 \ge t_0 $).

Нетрудно заметить, что средние значения
\[
\operatorname{ M } ( \varphi_0( t ) ) = - A^{ -1 } m_0,
\ \ \operatorname{ M } ( \varphi_{ 11 }( t ) ) = A^{ -2 } m_0 - A^{ -1 } m_1,
\ \ \operatorname{ M } ( \varphi_{ 12 }( t ) ) = A^{ -2 } m_0 - A^{ -1 } m_2.
\]
Поэтому формулу \eqref{NPI_DELC_exampleSolution} можно переписать в виде
$$
y = - A^{ -1 } m_0 +
\frac{1}{t} \left( A^{ -2 } m_0 - A^{ -1 } m_1 \right) +
\frac{1}{ t \ln t } \left( A^{ -2 } m_0 - A^{ -1 } m_2 \right) + r,
\ \ r = O \left( \frac{ 1 }{ t^2} \right) \ ( t \to + \infty ).
$$

\end{example}

Следующий пример иллюстрирует применение теоремы \ref{T:AsymptChart} об
асимптотическом характере формального частного решения (семейства решений)
\eqref{FSolution} основного вещественного векторного квазилинейного
обыкновенного дифференциального уравнения \eqref{DEV:main} в том случае,
когда матрица $ A $ не имеет кратных собственных значений, а класс
$ \mathcal K_3(A) $ состоит из РПП функций с конечными спектрами.

\begin{example}

Пусть числа $ t_0 > 1 $, $ n := m := 2 $, $ k_0 := 1 $, $ k := 3 $,
$ \varpi := 1,1 $, матрица $ A \in \mathbb R^{ 2 \times 2 } $ и выполнено
условие
\begin{equation*}
\inf_{ \substack{ \gamma \in \Gamma_2 \\ \lambda \in
\Lambda( A ) \cup \Delta( A ) } }
\left| \imath \gamma - \lambda \right| > 0,
\end{equation*}
где множество
\begin{equation*}
\Gamma_2 := \left\{ \left. \sum_{ s=1 }^4 k_s \omega_s \, \right|
k_s \in \mathbb Z, \ \omega_s \in \Omega_1
\ \left( s = \overline{ 1, 4 } \right) \right\},
\ \ \Omega_1 := \left\{ 1, \ \sqrt 2, \ e, \ \pi \right\}.
\end{equation*}
В качестве класса $ \mathcal K_3(A) $ возьмём множество вещественных РПП
функций обладающих конечными спектрами, являющимися подмножествами
$ \Gamma_2 $.
Зададим ранги степенно-логарифмических функций так, как указано в примере
\ref{E:rang}.

Рассмотрим в области $ D := \left\{ ( t, { y } ) \ \left| \ t \in I,
\ { y } \in \mathbb R^{ 2 \times 1 }, \right.
|| { y } - \varphi_0( t ) || \le a \right\} $
(число $ a \in \mathbb R^+ $) векторное дифференциальное уравнение
\begin{equation}
\frac{ dy }{ dt } = A y +
\begin{pmatrix}
  \sin t \\
  \cos e t \\
\end{pmatrix}
+ \frac{1}{t}
\begin{pmatrix}
  \tilde \gamma_1 ( y )_1 + ( y )_2^2 \cos \pi t \\
  \tilde \gamma_1 ( y )_2 \\
\end{pmatrix}
+ \frac{1}{ t \ln t }
\begin{pmatrix}
  \tilde \gamma_2 ( y )_1 \\
  \tilde \gamma_2 ( y )_2 + ( y )_1 ( y )_2 \sin \sqrt{2} \, t \\
\end{pmatrix},
                                                    \label{DEV:example}
\end{equation}
причём вектор--функция $ \varphi_0( t ) $ --- единственное частное решение из
класса $ \mathcal K_3^{ 2 \times 1 }(A) \cap C^1( I ) $ укороченного
векторного дифференциального уравнения соответствующего уравнению
\eqref{DEV:example}:
\begin{equation*}
\frac{ d \varphi_0( t ) }{ dt } = A \varphi_0( t ) +
\begin{pmatrix}
  \sin t \\
  \cos e t \\
\end{pmatrix},
\end{equation*}
числа $ \tilde \gamma_1, \ \tilde \gamma_2 \in \mathbb R $.

Очевидно, что в рассматриваемом случае
$$
\rho_s = s \ ( s \in \mathbb N ), \ \varkappa_1 = 2,
\ \nu_{ 11 }(t) \equiv \varepsilon_1(t) \equiv \frac{1}{t},
\ \nu_{ 12 }(t) \equiv \varepsilon_2(t) \equiv \frac{1}{ t \ln t }.
$$

Для нахождения вектор--функций $  \varphi_{ 11 }(t) $,
$  \varphi_{ 12 }(t) $ воспользуемся замечанием \ref{R:first_terms} в
случае класса $ \mathcal K_3(A) $.
Найдём их как единственные частные решения из класса
$ \mathcal K_3^{ 2 \times 1 }(A) \cap C^1( I ) $ векторных
дифференциальных уравнений из совокупности
\begin{equation*}
\left[
\begin{array}{l}
\dfrac{ d \varphi_{ 11 }(t) }{ dt } = A \varphi_{ 11 }(t) +
\begin{pmatrix}
  \tilde \gamma_1 ( \varphi_0( t ) )_1 +
  \left( \varphi_0( t ) \right)_2^2 \cos \pi t \\
  \tilde \gamma_1 ( \varphi_0( t ) )_2 \\
\end{pmatrix},
                                                                        \\
\dfrac{ d \varphi_{ 12 }(t) }{ dt } = A \varphi_{ 12 }(t) +
\begin{pmatrix}
  \tilde \gamma_2 ( \varphi_0( t ) )_1 \\
  \tilde \gamma_2 ( \varphi_0( t ) )_2 + \left( \varphi_0( t ) \right)_1
\left( \varphi_0( t ) \right)_2 \sin \sqrt{2} \, t \\
\end{pmatrix}.
\end{array}
\right.
\end{equation*}

И так, выполнены все условия теоремы \ref{T:AsymptChart},
следовательно, для достаточно большого числа $ t_0 $ у векторного
дифференциального уравнения \eqref{DEV:example} на промежутке $ I $
существует хотя бы одно частное решение вида
$$
y = \varphi_0( t ) + \frac{1}{t} { \varphi }_{ 11 }(t) +
\frac{1}{ t \ln t } { \varphi }_{ 12 }(t) + r,
\ \ r = O \left( \frac{ 1 }{ t^2} \right) \ ( t \to + \infty ).
$$
Более того, погрешность $ r $ зависит от стольких достаточно малых по
абсолютной величине произвольных скалярных постоянных, сколько имеется
индексов $ j \in \left\{ 1, 2 \right\} $ таких, что справедливы условия
\begin{equation}
\lim_{ t \to + \infty } \Re \, \left( \lambda_j(A) t +
\left( \gamma_{ j 10 } + 2,9 \right) \ln t +
\gamma_{ j 01 } \ln \ln t \right) \ne + \infty,
                                            \label{E:C:No_of_parameters_O}
\end{equation}
где числа $ \gamma_{ j 10 } $, $ \gamma_{ j 01 } \in \mathbb C $.
Воспользовавшись формулами \eqref{gamma_jp, |p|=1} для нахождения этих чисел,
получим
\begin{gather*}
\gamma_{ j 10 } = \operatorname{ M } \left( \left( P_0^{ -1 }
\begin{pmatrix}
  \tilde \gamma_1 & 2 \left( \varphi_0( t ) \right)_2 \cos \pi t \\
  0 & \tilde \gamma_1 \\
\end{pmatrix}
P_0 \right)_{ jj } \right),
\\
\gamma_{ j 01 } = \operatorname{ M } \left( \left( P_0^{ -1 }
\begin{pmatrix}
  \tilde \gamma_2 & 0 \\
  \left( \varphi_0( t ) \right)_2 \sin \sqrt{2} \, t &
  \tilde \gamma_2 + \left( \varphi_0( t ) \right)_1 \sin \sqrt{2} \, t \\
\end{pmatrix}
P_0 \right)_{ jj } \right)
\ \left( j = \overline{ 1, 2 } \right),
\end{gather*}
где матрица $ P_0 \in \mathbb C^{ 2 \times 2 } $ составлена из собственных
векторов матрицы $ A $.
Очевидно, что средние матрицы из этих выражений можно представить в виде
суммы постоянной диагональной и переменной матриц.
Выполнив действия с первым слагаемым, получим
\begin{gather*}
\gamma_{ j 10 } = \tilde \gamma_1 + \operatorname{ M } \left( \left( P_0^{ -1 }
\begin{pmatrix}
  0 & 2 \left( \varphi_0( t ) \right)_2 \cos \pi t \\
  0 & 0 \\
\end{pmatrix}
P_0 \right)_{ jj } \right),
\\
\gamma_{ j 01 } = \tilde \gamma_2 + \operatorname{ M } \left( \left( P_0^{ -1 }
\begin{pmatrix}
  0 & 0 \\
  \left( \varphi_0( t ) \right)_2 \sin \sqrt{2} \, t &
  \left( \varphi_0( t ) \right)_1 \sin \sqrt{2} \, t \\
\end{pmatrix}
P_0 \right)_{ jj } \right)
\ \left( j = \overline{ 1, 2 } \right).
\end{gather*}
Заметим, что вектор--функция
$ \varphi_0( t ) \in \mathcal K_3^{ 2 \times 1 }(A) $
является линейной комбинацией функций $ \cos \gamma t $, $ \sin \gamma t $
($ \gamma = 1, e $) с векторными коэффициентами из
$ \mathbb R^{ 2 \times 1 } $.
Поэтому у произведений её компонент на функции $ \cos \gamma t $,
$ \sin \gamma t $ ($ \gamma = \sqrt 2, \pi $) не будет средних значений.
Следовательно, у последних произведений матриц тоже не будет средних
значений.
Значит, числа $ \gamma_{ j 10 } = \tilde \gamma_1 $,
$ \gamma_{ j 01 } = \tilde \gamma_2 $
$ \left( j = \overline{ 1, 2 } \right) $.
В итоге условия \eqref{E:C:No_of_parameters_O} можно переписать в следующем
виде
\begin{equation*}
\lim_{ t \to + \infty } \left( \Re \, \lambda_j(A) t +
\left( \tilde \gamma_1 + 2,9 \right) \ln t +
\tilde \gamma_2 \ln \ln t \right) \ne + \infty
\ \left( j \in \left\{ 1, 2 \right\} \right).
\end{equation*}

\end{example}

Эта работа выполнена по теме имеющей номер в государственной регистрации
Украины 0109U003444.

Перспективы дальнейших изысканий по этой теме представляются в изучении
особых случаев.
(Например, случай кратных чисто мнимых собственных значений матрицы $ A $.)

\bibliography{Amelkin_KV,Other}

\begin{thebibliography}{10}
\def\selectlanguageifdefined#1{
\expandafter\ifx\csname date#1\endcsname\relax
\else\language\csname l@#1\endcsname\fi}
\providecommand*{\href}[2]{{\small #2}}
\providecommand*{\url}[1]{{\small #1}}
\providecommand*{\BibUrl}[1]{\url{#1}}
\providecommand{\BibAnnote}[1]{}
\providecommand*{\BibEmph}[1]{#1}
\providecommand*{\cyrdash}{\hbox to.8em{--\hss--}}
\providecommand*{\BibDash}{\ifdim\lastskip>0pt\unskip\nobreak\hskip.2em\fi
\cyrdash\hskip.2em\ignorespaces}

\bibitem{Bib:Kamenec-Podolsky96}
\selectlanguageifdefined{russian}
\BibEmph{Амелькин~К.~В., Костин~А.~В.} О методе малого переменного параметра
  исследования одного нелинейного дифференциального уравнения второго
  порядка~// Нелінійні крайові задачі математичної фізики та їх застосування:
  Збірник праць. \BibDash
\newblock Киев~: Институт математики НАН Украины, 1996. \BibDash
\newblock {\cyr\CYRS.}~14--15.

\bibitem{Bib:Mastersthesis}
\selectlanguageifdefined{russian}
\BibEmph{Амелькин~К.~В.} О методе малого переменного параметра в теории
  колебаний. \BibDash
\newblock 1997. \BibDash
\newblock Одесский Национальный университет им. И. И. Мечникова.

\bibitem{Bib:Chernovci98}
\selectlanguageifdefined{russian}
\BibEmph{Костин~А.~В., Амелькин~К.~В.} О методе малого переменного параметра в
  теории колебаний~// СУЧАСНІ ПРОБЛЕМИ МАТЕМАТИКИ : Матеріали Міжнародної
  наукової конференції. Частина 2. \BibDash
\newblock Киев~: Институт математики НАН Украины, 1998. \BibDash
\newblock {\cyr\CYRS.}~13--16.

\bibitem{Bib:Odessa00}
\selectlanguageifdefined{russian}
\BibEmph{Амелькин~К.~В., Костин~А.~В.} О колебаниях в квазилинейных системах
  дифференциальных уравнений с исчезающими при $ t \to + \infty $ параметрами
  при нелинейностях~// Диференціальні та інтегральні рівняння. Тези доповідей
  Міжнародної конференції, 12 – 14 вересня 2000. \BibDash
\newblock Одесса~: Астропринт, 2000. \BibDash
\newblock {\cyr\CYRS.}~9--10.

\bibitem{Bib:Chernovci01}
\selectlanguageifdefined{russian}
\BibEmph{Амелькин~К.~В., Костин~А.~В.} О методе малых переменных параметров~//
  Диференціальні рівняння і нелінійні коливання. Тези доповідей Міжнародної
  конференції, 27–29 серпня 2001. \BibDash
\newblock Киев~: Институт математики НАН Украины, 2001. \BibDash
\newblock {\cyr\CYRS.}~8--9.

\bibitem{Bib:Chernovci02}
\selectlanguageifdefined{russian}
\BibEmph{Костин~А.~В., Амелькин~К.~В.} О методе малых переменных параметров для
  двумерных квазилинейных дифференциальных систем~// Диференціальні рівняння і
  нелінійні коливання. Секція 3: Праці Українського математичного конгресу –
  2001. \BibDash
\newblock Киев~: Институт математики НАН Украины, 2002. \BibDash
\newblock {\cyr\CYRS.}~34--42.

\bibitem{Bib:Odessa01}
\selectlanguageifdefined{russian}
\BibEmph{Амелькин~К.~В., Костин~А.~В.} О расщеплении линейных однородных систем
  обыкновенных дифференциальных уравнений~// \BibEmph{Вісник Одеського
  державного університету. Фіз.-мат. науки}. \BibDash
\newblock 2001. \BibDash
\newblock \CYRT.~6, {\cyr\textnumero}~3. \BibDash
\newblock {\cyr\CYRS.}~36--43. \BibDash
\newblock URL: \BibUrl{http://www.scribd.com/embeds/114123650/
  content?start_page=1&view_mode=scroll&access_key=key-n5giqabsg2jisgp3pum}.

\bibitem{Bib:RNASU04}
\selectlanguageifdefined{russian}
\BibEmph{Костин~А.~В., Амелькин~К.~В.} О методе малых переменных параметров для
  квазилинейных дифференциальных систем~// \BibEmph{Доклады НАН Украины}.
  \BibDash
\newblock 2004. \BibDash
\newblock \CYRT.~6. \BibDash
\newblock {\cyr\CYRS.}~13--18.

\bibitem{Bib:Kostin_DE65}
\selectlanguageifdefined{russian}
\BibEmph{Костин~А.~В.} К вопросу о существовании у системы обыкновенных
  дифференциальных уравнений ограниченных частных решений и частных решений,
  стремящихся к нулю при $ t \to + \infty $~// \BibEmph{Дифференциальные
  уравнения}. \BibDash
\newblock 1965. \BibDash
\newblock \CYRT.~1, {\cyr\textnumero}~5. \BibDash
\newblock {\cyr\CYRS.}~585--604.

\bibitem{Bib:Kostin_DE67_I}
\selectlanguageifdefined{russian}
\BibEmph{Костин~А.~В.} Об асимптотических рядах в теории нелинейных систем
  обыкновенных дифференциальных уравнений. 1~// \BibEmph{Дифференциальные
  уравнения}. \BibDash
\newblock 1967. \BibDash
\newblock \CYRT.~3, {\cyr\textnumero}~6. \BibDash
\newblock {\cyr\CYRS.}~875–889.

\bibitem{Bib:Kostin_DE67_II}
\selectlanguageifdefined{russian}
\BibEmph{Костин~А.~В.} Об асимптотических рядах в теории нелинейных систем
  обыкновенных дифференциальных уравнений. 2~// \BibEmph{Дифференциальные
  уравнения}. \BibDash
\newblock 1967. \BibDash
\newblock \CYRT.~3, {\cyr\textnumero}~7. \BibDash
\newblock {\cyr\CYRS.}~1070–1077.

\bibitem{Bib:Kostin_PHDT}
\selectlanguageifdefined{russian}
\BibEmph{Костин~А.~В.} Устойчивость и асимптотика квазилинейных неавтономных
  дифференциальных систем. \BibDash
\newblock Одесса~: ОГУ, 1984. \BibDash
\newblock {\cyr\CYRS.}~94.

\bibitem{Bib:Kostin_DE87}
\selectlanguageifdefined{russian}
\BibEmph{Костин~А.~В.} Асимптотические разложения исчезающих и ограниченных
  решений обыкновенных дифференциальных уравнений~// \BibEmph{Дифференциальные
  уравнения}. \BibDash
\newblock 1987. \BibDash
\newblock \CYRT.~23, {\cyr\textnumero}~7. \BibDash
\newblock {\cyr\CYRS.}~1629–1632.

\bibitem{Bib:Kostin_DNANU95}
\selectlanguageifdefined{russian}
\BibEmph{Костин~А.~В., Кореновский~А.~А.} Асимптотические разложения
  обобщённого типа $ o $\nobreakdash-решений квазилинейной неавтономной системы
  дифференциальных уравнений~// \BibEmph{Доклады НАН Украины}. \BibDash
\newblock 1995. \BibDash
\newblock \CYRT.~10. \BibDash
\newblock {\cyr\CYRS.}~13–15.

\bibitem{Bib:Abgaryan2008}
\selectlanguageifdefined{russian}
\BibEmph{Абгарян~К.~А.} Матричное исчисление с приложениями в теории
  динамических систем. Учебное пособие для вузов. \BibDash
\newblock 2-е {\cyr\cyri\cyrz\cyrd.} \BibDash
\newblock Москва~: Вузовская книга, 2008. \BibDash
\newblock {\cyr\CYRS.}~544.

\bibitem{Bib:Mitropolskiy}
\selectlanguageifdefined{russian}
\BibEmph{Митропольский~Ю.~А.} Нелинейная механика. Асимптотические методы.
  \BibDash
\newblock Киев~: Институт математики НАН Украины, 1995. \BibDash
\newblock {\cyr\CYRS.}~396.

\bibitem{Bib:Vazov}
\selectlanguageifdefined{russian}
\BibEmph{Вазов~В.} Асимптотические разложения решений обыкновенных
  дифференциальных уравнений. \BibDash
\newblock Москва~: Мир, 1968. \BibDash
\newblock {\cyr\CYRS.}~464.

\bibitem{Bib:Demidovich}
\selectlanguageifdefined{russian}
\BibEmph{Демидович~Б.~П.} Лекции по математической теории устойчивости.
  \BibDash
\newblock Москва~: Наука, 1967. \BibDash
\newblock {\cyr\CYRS.}~472.

\bibitem{Bib:Malkin}
\selectlanguageifdefined{russian}
\BibEmph{Малкин~И.~Г.} Некоторые задачи теории нелинейных колебаний. \BibDash
\newblock Москва~: Гостехиздат, 1956. \BibDash
\newblock {\cyr\CYRS.}~491.

\bibitem{Bib:Levitan}
\selectlanguageifdefined{russian}
\BibEmph{Левитан~Б.~М.} Почти--периодические функции. \BibDash
\newblock Москва~: ГИТТЛ, 1953. \BibDash
\newblock {\cyr\CYRS.}~396.

\bibitem{Bib:Sokolov}
\selectlanguageifdefined{russian}
\BibEmph{Соколов~Ю.~Д.} Об асимптотических решениях дифференциальных
  уравнений~// \BibEmph{Сборник трудов Киевского инженерно-строительного
  института}. \BibDash
\newblock 1948. \BibDash
\newblock {\cyr\textnumero}~8. \BibDash
\newblock {\cyr\CYRS.}~62--74.

\bibitem{Bib:Persidskiy}
\selectlanguageifdefined{russian}
\BibEmph{Персидский~К.~П.} О характеристических числах дифференциальных
  уравнений~// \BibEmph{Известия АН Казахской ССР. Серия математики и
  механики}. \BibDash
\newblock 1947. \BibDash
\newblock \CYRT.~42, {\cyr\textnumero}~1. \BibDash
\newblock {\cyr\CYRS.}~5--47.

\bibitem{Bib:Cotton}
\selectlanguageifdefined{english}
\BibEmph{Cotton~E.} Sur les solutions asymptotiques des equations
  differentielles~// \BibEmph{Annales scientifiques de l'Ecole Normale
  Superieure}. \BibDash
\newblock 1911. \BibDash
\newblock Vol.~28, no.~3. \BibDash
\newblock P.~473--521.

\end{thebibliography}

\end{document}